\documentclass[3p,12pt]{elsarticle}

\usepackage{fullpage}
\usepackage{amsmath,amsfonts,amssymb}
\usepackage{amsthm}
\usepackage{newtxtext,newtxmath}
\usepackage[toc]{appendix}
\usepackage[T1]{fontenc}
\usepackage{textcomp}
\usepackage{graphicx,epstopdf}
\usepackage[position=top,labelfont=normalfont,textfont=normalfont,singlelinecheck=off,justification=raggedright]{subfig}
\usepackage{floatrow}
\usepackage[dvipsnames]{xcolor}
\usepackage{setspace}
\usepackage{enumerate,etaremune}
\usepackage{tabularx,multirow}
\usepackage{framed}
\usepackage{hhline}
\usepackage[pdftex,unicode,bookmarks=false]{hyperref}
\usepackage{url}
\hypersetup{
  colorlinks,
  citecolor=Blue, linkcolor=Blue, urlcolor=Blue
}

\usepackage{tikz}
\usetikzlibrary{shapes.geometric}
\usetikzlibrary{shapes,arrows}
\usetikzlibrary{plotmarks}

\usetikzlibrary{calc}

\usepackage{algorithm}
\usepackage{algorithmicx,algpseudocode}
\usepackage{cases}

\newcommand{\eg}{{e.g.}}
\newcommand{\ie}{{i.e.}}

\newcommand{\etal}{{\it et al.}}
\newcommand{\tensor}[1]{\ensuremath{\boldsymbol{#1}}}

\newcommand{\eps}{\varepsilon}
\newcommand{\od}{\mathrm{d}}
\newcommand{\pd}{\partial}
\newcommand{\Del}{\mathrm{\Delta}}

\newcommand{\cn}{\mathrm{N}}

\DeclareMathOperator{\grad}{\nabla}
\DeclareMathOperator{\diver}{\grad\cdot}
\DeclareMathOperator{\symgrad}{\nabla^{s}}

\DeclareMathOperator{\dyad}{\otimes}
\DeclareMathOperator{\sign}{sign}
\newsavebox{\dotbox}

\theoremstyle{remark}
\newtheorem{remark}{Remark}
\newtheorem*{remark*}{Remark}

\newcommand{\revised}[1]{{\color{black} #1}}

\newcolumntype{L}[1]{>{\raggedright\let\newline\\arraybackslash\hspace{0pt}}m{#1}}
\newcolumntype{C}[1]{>{\centering\let\newline\\arraybackslash\hspace{0pt}}m{#1}}
\newcolumntype{R}[1]{>{\raggedleft\let\newline\\arraybackslash\hspace{0pt}}m{#1}}

\allowdisplaybreaks

\biboptions{sort&compress,square,comma,numbers}
\AtBeginDocument{\hypersetup{citecolor=MidnightBlue,linkcolor=MidnightBlue,urlcolor=MidnightBlue}}
\let\today\relax

\begin{document}

\begin{frontmatter}

\title{A phase-field method for modeling cracks with frictional contact}

\author[HKU]{Fan Fei}
\author[HKU]{Jinhyun Choo\corref{corr}}
\ead{jchoo@hku.hk}

\cortext[corr]{Corresponding Author}

\address[HKU]{Department of Civil Engineering, The University of Hong Kong, Pokfulam, Hong Kong}

\journal{$\,$}

\begin{abstract}
We introduce a phase-field method for continuous modeling of cracks with frictional contacts.
Compared with standard discrete methods for frictional contacts, the phase-field method has two attractive features: (1) it can represent arbitrary crack geometry without an explicit function or basis enrichment, and (2) it does not require an algorithm for imposing contact constraints.
The first feature, which is common in phase-field models of fracture, is attained by regularizing a sharp interface geometry using a surface density functional.
The second feature, which is a unique advantage for contact problems, is achieved by a new approach that calculates the stress tensor in the regularized interface region depending on the contact condition of the interface.
Particularly, under a slip condition, this approach updates stress components in the slip direction using a standard contact constitutive law, while making other stress components compatible with stress in the bulk region to ensure non-penetrating deformation in other directions.
We verify the proposed phase-field method using stationary interface problems simulated by discrete methods in the literature.
Subsequently, by allowing the phase field to evolve according to brittle fracture theory, we demonstrate the proposed method's capability for modeling crack growth with frictional contact.
\end{abstract}


\end{frontmatter}

\section{Introduction}
Frictional cracks are ubiquitous in natural and manufactured systems.
For example, in the Earth's crust, frictional cracks appear over a wide range of scales \revised{from the micrometer scale~\cite{Vajdova2010,Tjioe2015,Tjioe2016} to the millimeter scale~\cite{Wong2009a,Wong2009b,White2014} to the kilometer scale~\cite{Sanz2007,Liu2009,Liu2013}}.
They are also widespread in many branches of science and engineering including material sciences and civil and mechanical engineering.
Accordingly, the numerical modeling of motion and friction in crack surfaces has long been an important subject,
and there is a large body of literature on computational contact mechanics (see, for example, Laursen~\cite{Laursen2003}, Wriggers~\cite{Wriggers2006}, and references therein).

At present, standard numerical methods treat frictional cracks as discrete discontinuities subjected to constraints on contact conditions.
The discontinuities should be aligned with element boundaries in classical finite element methods, while they can be embedded inside elements in modern methods \revised{such as the assumed enhanced strain (AES) method~\cite{Simo1993,Foster2007,Borja2008} and the extended finite element method (XFEM)~\cite{Dolbow2001,Liu2008,Sanborn2011}.}
Irrespective of their alignment with elements, the contact surfaces must satisfy a set of constraints including the no-penetration constraint under compression.
Imposing these constraints on discrete interfaces, however, is an outstanding challenge in computational contact mechanics.
A large number of studies have addressed this challenge by employing various types of algorithms such as \revised{the Lagrange multiplier method~\cite{Simo1985,Carpenter1991,Liu2010b}, the penalty method~\cite{Peric1992,Khoei2007,Liu2010a,Mueller-Hoeppe2012}, the augmented Lagrangian method \cite{Simo1992,Elguedj2007,Bussetta2012}, and the Nitsche method~\cite{Dolbow2009,Annavarapu2012,Annavarapu2013,Annavarapu2014}.}
Nevertheless, the optimal way to treat these contact constraints is yet an unresolved issue.
Also importantly, numerical methods for frictional contacts require significant effort for implementation, especially when one wants to accommodate complex interface geometry.
For these reasons, a simple numerical method is desired that can model frictional contact problems with low implementation cost.

In this paper, we propose a phase-field method for frictional crack problems that can efficiently handle complex crack geometry and contact conditions.
Phase-field modeling is a continuous (as opposed to discrete) approach to interface problems that approximates a sharp interface as a region where the phase field attains a certain value.
Its upshot is that one can represent an interface without any function or algorithm for describing its geometry, which is highly advantageous when the geometry is complex and may evolve with time.
For this reason, phase-field modeling has found widespread applications in a variety of scientific and engineering problems.
Recently, it has enjoyed considerable success in computational modeling of a wide range of fracturing processes \revised{under mechanical~\cite{Miehe2010,Borden2012,Choo2018a,Geelen2019}, thermal~\cite{Bourdin2014,Miehe2015,Na2018}, hydraulic~\cite{Lee2016b,Santillan2017,Ha2018}, and chemical loads~\cite{Miehe2016,Zhang2016,Choo2018b}.}
Nevertheless, to our knowledge, the present work is the first attempt to apply a phase-field approach to cracks with frictional contact.

The key idea of the proposed phase-field method is to incorporate contact behaviors and constraints through suitable calculation of the stress tensor in the regularized interface region.
In existing phase-field models of fracture, the stress tensor in the interface region is either degraded or maintained according to the sign of some part of the strain tensor.
This way roughly considers a no-contact condition and a stick contact condition, but these two conditions are usually not distinguished in a manner consistent with contact mechanics.
Furthermore, and perhaps more importantly, phase-field fracture models have not incorporated a slip contact condition in which relative motion between interacting surfaces takes places according to friction.
Note that a slip condition is a major challenge in computational contact mechanics because it requires one to model slip behavior while imposing the no-penetration contact in non-slip directions.
In this work, we propose a new approach that incorporates and distinguishes between all contact conditions by a proper calculation of the stress tensor in the interface region.
Our use of stress tensor is consistent with other types of smeared crack formulations for frictional cracks (\eg~Borja~\cite{Borja2000}), but our way to calculate stress is completely different from existing smeared methods as it builds on a stress calculation procedure in phase-field modeling of fracture.

Compared with standard discrete methods for frictional contacts, the proposed phase-field method has two attractive features.
First, it can represent arbitrary interface geometry without an explicit function or enriched basis functions, which is indeed a hallmark of all phase-field methods.
Second, it can accommodate contact constraints without a dedicated algorithm, which is a unique advantage for contact problems.
This new feature is attained by making the components of the stress tensor in the non-slip directions compatible between the interface and bulk regions.
Notably, this way is a modification of the volumetric--deviatoric decomposition approach proposed by Amor \etal~\cite{Amor2009} for considering the no-penetration constraint in unilateral frictionless contact, \revised{as well as the directional decomposition approach proposed by Steinke and Kaliske~\cite{Steinke2019} for distinguishing between crack normal and tangential directions.}
The rest components of the stress tensor, which are relevant to stick/slip behavior, are determined using a standard constitutive law for frictional cracks.
In this way, the proposed phase-field method translates a discrete problem with constraints into a continuous problem with multiple constitutive responses.
As the continuous problem can be solved by the standard finite element method,
it would require significantly less effort for implementation as compared with other methods that can simulate frictional cracks passing through the interior of elements.

The paper is organized as follows.
In Section~\ref{sec:formulation}, we develop a phase-field formulation for a boundary-value problem that involves frictional contact.
Following a standard diffuse approximation of crack geometry in phase-field modeling of fracture,
we introduce a new approach that explicitly considers and calculates the stress tensor in the interface region according to the contact condition of the interface.
In Section~\ref{sec:discretization}, we discretize the phase-field contact formulation using the standard finite element method.
In doing so, we present algorithms for calculating the stress, stress--strain tangent, and unit normal and slip vectors at quadrature points in the regularized interface region.
In Section~\ref{sec:examples}, we verify the proposed phase-field method using stationary interface problems that have been simulated by discrete methods in the literature.
We then combine the proposed method with an evolution equation for brittle fracture, and demonstrate the method's capability for modeling crack growth with frictional contact.
In Section~\ref{sec:closure}, we conclude the work.

\section{Phase-field formulation for cracks with frictional contact}
\label{sec:formulation}
In this section, we develop a phase-field formulation for continuous modeling of frictional cracks in solids.
The formulation builds on methods originally developed for phase-field modeling of crack propagation,
but it can be useful for modeling general frictional interfaces in solids.
For this reason, we will focus on the use of a phase-field approach to geometric approximation of frictional contact problems, without delving into aspects of fracture mechanics.
Also, to make the following presentation simple, we will assume without loss of generality that the material is free of inertial and body forces, isotropic, elastic, and geometrically and materially linear.
If necessary, these assumptions may be relaxed in standard ways in solid mechanics.

\subsection{Problem statement and phase-field approximation}
Consider a domain $\Omega\in\mathbb{R}^{\mathrm{dim}}$ in a $\mathrm{dim}$-dimensional space with external boundary $\partial\Omega$.
The boundary is partitioned into the displacement (Dirichlet) boundary, $\partial_{u}\Omega$, and the traction (Neumann) boundary, $\partial_{t}\Omega$, such that $\overline{\partial_{u}\Omega\cap\partial_{t}\Omega}=\emptyset$ and $\overline{\partial_{u}\Omega\cup\partial_{t}\Omega}=\partial\Omega$. The domain has a set of internal discontinuities, which is denoted by $\Gamma$.
A discontinuity has two surfaces that are either separated or in contact.
When the two surfaces are in contact, relative motion may or may not exist between them depending on the magnitudes of tractions and frictions therein.

\begin{figure}[htbp]
  \centering
  \includegraphics[width=0.9\textwidth]{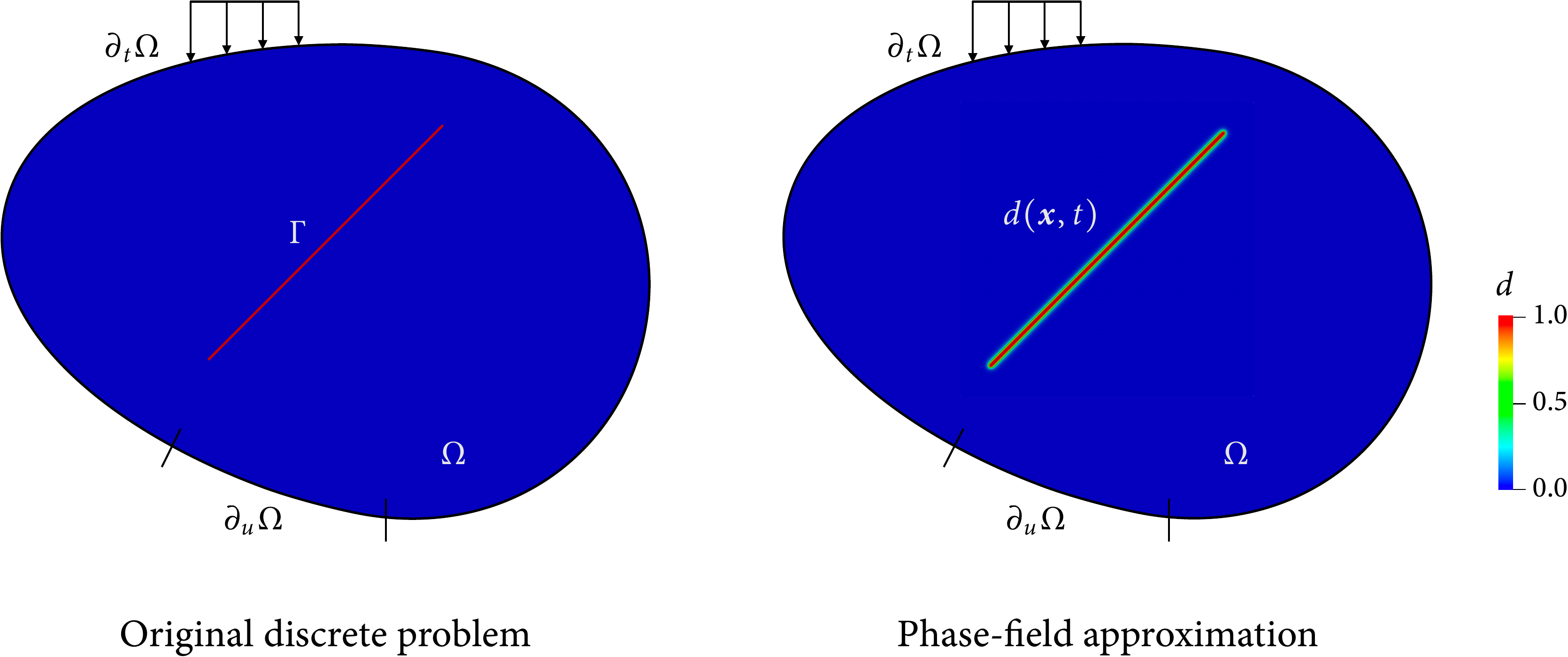}
  \caption{Phase-field approximation of a discrete problem with a frictional crack.
  The sharp discontinuity $\Gamma$ in the left figure is diffusely approximated by the phase-field variable $d$ as illustrated in the right figure.}
  \label{fig:phase-field-approximation}
\end{figure}

We begin our formulation by approximating the discontinuities' geometry using a standard approach in phase-field modeling of fracture.
Let us define the phase-field variable, $d\in[0,1]$, such that it denotes a fully discontinuous (interface) region by $d=1$ and a fully continuous (bulk) region by $d=0$.
We then introduce a surface density functional for length regularization of the sharp geometry of $\Gamma$.
Among several forms of the functional proposed in the literature, here we adopt the most popular one, given by
\begin{align}
  \gamma(d,\grad d) := \frac{1}{2}\left(\frac{d^{2}}{L} + L\grad{d}\cdot\grad{d} \right)\,.
\end{align}
Here, $L$ is the length parameter introduced for regularization of sharp geometry, which determines the size of the diffuse approximation zone.
Figure~\ref{fig:phase-field-approximation} illustrates how this phase-field approach approximates the original domain with sharp discontinuity.
Note that the diffuse approximation naturally gives rise to regions in which $0<d<1$.

Once discontinuities have been diffusely approximated in the way described above,
a continuous version of the problem can be stated as follows.
Find the displacement field in this domain, $\tensor{u}$, that satisfies the balance of linear momentum
\begin{align}
  \diver\,\tensor{\sigma}(\tensor{\eps}) = \tensor{0} \quad\text{in}\;\;\Omega\,,
\end{align}
where $\tensor{\sigma}$ is the (Cauchy) stress tensor
and $\tensor{\eps}$ is the infinitesimal strain tensor defined as the symmetric gradient of $\tensor{u}$.
The boundary conditions of this problem are given by
\begin{align}
  \tensor{u} = \hat{\tensor{u}} \quad&\text{on}\;\;\partial_{u}\Omega\,,\\
  \tensor{\upsilon}\cdot\tensor{\sigma} = \hat{\tensor{t}} \quad&\text{on}\;\;\partial_{t}\Omega\,,
\end{align}
where $\hat{\tensor{u}}$ and $\hat{\tensor{t}}$ are prescribed boundary conditions of displacement and traction vectors, respectively, and $\tensor{\upsilon}$ is the outward unit normal vector at the domain boundary.
Note that no boundary condition is imposed on $\Gamma$ as the discontinuities have already been smeared in the domain $\Omega$ through the phase-field approximation described above.

As the whole domain is now regarded as a continuum, stress tensors in the bulk and interface systems should be continuously interpolated.
For this interpolation, we introduce a function of the phase field, $g(d)$, that satisfies
\begin{align}
  g(d) \in[0,1]\,, \quad
  g(0) = 1\,, \quad
  g(1) = 0\,, \quad
  g'(d)<0\,.
\end{align}
This function is usually called the degradation function in phase-field modeling of fracture.
In this work, we use the most common form of $g(d)$ in the literature, given by
\begin{align}
  g(d) = (1-d)^{2}\,,
\end{align}
Using this function, we can express the stress tensor in the domain as
\begin{align}
  \tensor{\sigma} = g(d)\tensor{\sigma}_{\text{bulk}} + [1 - g(d)]\tensor{\sigma}_{\text{interface}}\,.
  \label{eq:stress-interpolation}
\end{align}
where $\tensor{\sigma}_{\text{bulk}}$ and $\tensor{\sigma}_{\text{interface}}$ are stress tensors
in the bulk and interface systems, respectively.
It is noted that Eq.~\eqref{eq:stress-interpolation} is a generalization of the way in which stress is calculated in phase-field models of fracture.

In addition, the phase-field variable is postulated to satisfy the following partial differential equation
\begin{align}
  g'(d)\mathcal{H} + G_{c}\left(\frac{d}{L} - L\diver\grad{d}\right) = 0 \quad\text{in}\;\;\Omega\,,
  \label{eq:phase-field}
\end{align}
which is also adopted from phase-field modeling of fracture.
Here, $\mathcal{H}$ and $G_{c}$ are positive parameters corresponding to the crack driving force and the critical fracture energy in the context of phase-field modeling of fracture.
As such, for simulation of growing cracks, they may be calculated according to a phase-field formulation for fracture.
However, if the crack interface is assumed to be stationary, one may take any positive values for these two parameters for initialization of the phase field.
In this case, Eq.~\eqref{eq:phase-field} is solved only once in the beginning of the problem to initialize the phase-field variable, and the phase field remains constant throughout the course of loading.

\subsection{Calculation of stress tensors according to contact conditions}
\label{subsec:stress-calculation}
So far, the only difference between our formulation and the most standard phase-field formulation for fracture is that here we have explicitly considered the stress tensor in the interface system, $\tensor{\sigma}_{\text{interface}}$.
This modification has been made to incorporate contact-dependent mechanical responses of the interface system into the phase-field formulation.
In the following, we propose a specific procedure for calculating the bulk and interface stress tensors in Eq.~\eqref{eq:stress-interpolation} according to the contact condition of the interface.

First, we calculate the bulk stress through a standard stress--strain relationship in continuum mechanics, namely
\begin{align}
  \tensor{\sigma}_{\text{bulk}} = \mathbb{C}_{\text{bulk}}:\tensor{\eps}\,,
\end{align}
where $\mathbb{C}_{\text{bulk}}$ is the fourth-order stress--strain tangent tensor of the bulk region.
For a linear elastic material, the stress--strain tangent is given by
\begin{align}
  \mathbb{C}_{\text{bulk}} = \mathbb{C}^{\rm e} := \lambda \tensor{1}\dyad\tensor{1} + 2G\mathbb{I}\,.
\end{align}
Here, $\lambda$ and $G$ are the Lam\'{e} parameters which can be converted into Young's modulus $E$ and Poisson's ratio $\nu$, $\tensor{1}$ is the second-order identity tensor,
and $\mathbb{I}$ is the fourth-order symmetric identity tensor.
In short, the bulk stress tensor is evaluated as usual.

\begin{figure}[htbp]
  \centering
  \includegraphics[width=0.9\textwidth]{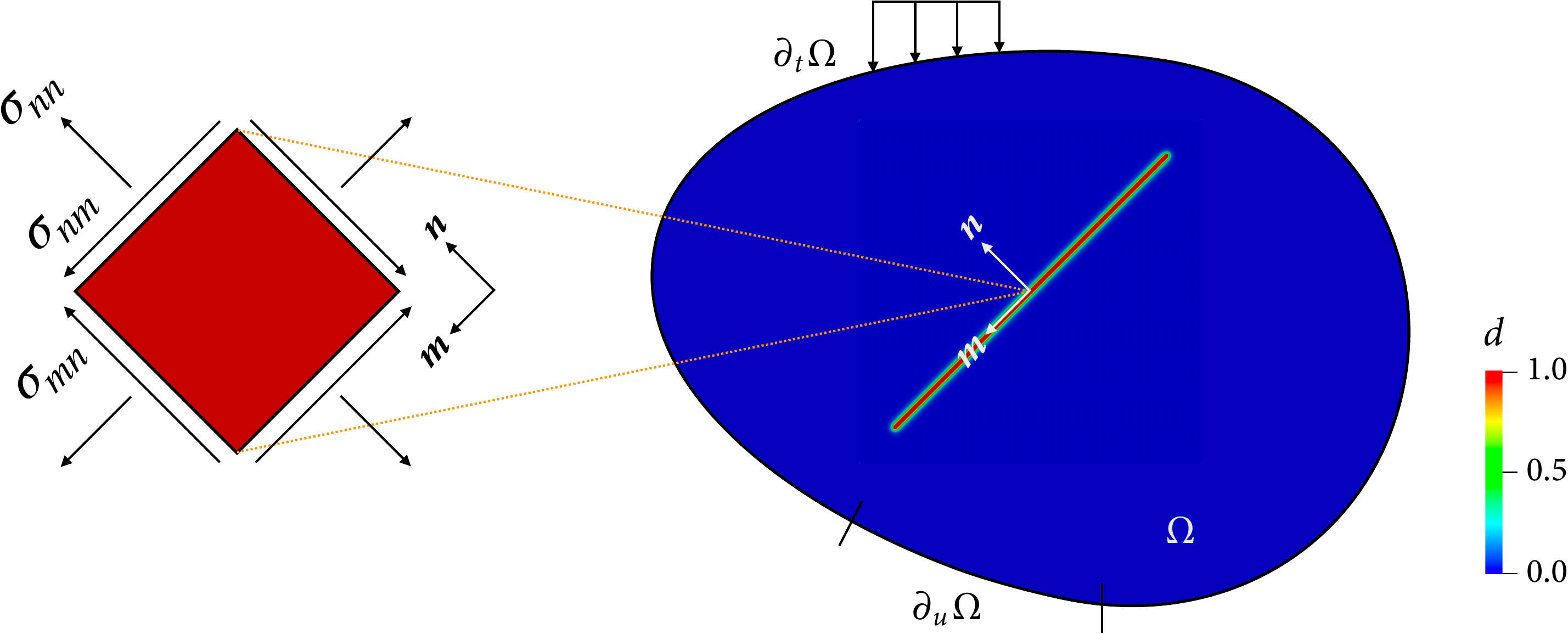}
  \caption{Definition of an interface-oriented coordinate system in a 2D domain.
  Vectors $\tensor{n}$ and $\tensor{m}$ denote unit vectors in the interface normal and tangential/slip directions, respectively. Definitions of some stress components in the interface-oriented coordinate system are also shown.}
  \label{fig:interface-coord}
\end{figure}

Next, we propose a new way to calculate the interface stress tensor depending on the contact condition of the interface system.
To identify the contact condition, we introduce a coordinate system that is oriented with respect to the interface normal and tangential directions.
Figure~\ref{fig:interface-coord} depicts this interface-oriented coordinate system in a 2D domain.
Hereafter, we denote by $\tensor{n}$ the unit vector in the interface normal direction and by $\tensor{m}$ the unit vector in the slip direction.
These unit vectors are assumed to be known for now.

In the interface-oriented coordinate system, the normal strain along the interface normal direction is calculated as
\begin{align}
  \eps_{\cn} &\equiv \eps_{nn} = \tensor{\eps}:(\tensor{n}\dyad\tensor{n})\,.
  \label{eq:interface-normal-strain}
\end{align}
This strain can be used to distinguish between contact and no-contact conditions.
By definition, this strain plays the role of the gap function in classical contact mechanics.

When $\eps_{\cn} > 0$, the interface has a gap between its two surfaces, which corresponds to an open (non-contacting) crack in phase-field modeling of fracture.
In this case, the interface system is stress-free, \ie~$\tensor{\sigma}_{\text{interface}}=\tensor{0}$.
We thus evaluate the stress tensor under this no-contact condition as
\begin{align}
  \tensor{\sigma} = g(d)\tensor{\sigma}_{\text{bulk}}\,.
\end{align}
It is noted that the above expression is the same as the stress equation for an open crack in phase-field models of fracture.

By contrast, when $\eps_{\cn} \leq 0$, the interface is considered being in contact.
The contact condition of a cohesive--frictional interface is either a stick condition or a slip condition.
The distinction between stick and slip conditions can be made by the following yield function
\begin{align}
  f := |\tau| - \mu p_{\cn} \leq 0\,.
  \label{eq:slip-yield-function}
\end{align}
Here, $\mu$ is the friction coefficient of the interface,
and $p_{\cn}$ is the contact normal pressure, defined as $p_{\cn} \equiv -\sigma_{nn} = -\tensor{\sigma}:(\tensor{n}\dyad\tensor{n})$ in accordance to Eq.~\eqref{eq:interface-normal-strain}.
Note that $p_{\cn} \geq 0$ whenever $\eps_{\cn} \leq 0$.
Lastly, $\tau$ is the resolved shear stress in the interface, which can be calculated in this continuum formulation as (see Fig.~\ref{fig:interface-coord})
\begin{align}
  \tau \equiv \sigma_{nm} = \tensor{\sigma}:(\tensor{n}\dyad\tensor{m})\,.
\end{align}
Note that $\sigma_{nm}=\sigma_{mn}$ by the symmetry of the stress tensor.
A contacting interface is under a stick condition when $f<0$, whereas it is under a slip condition when $f=0$.
Therefore, stick and slip conditions are distinguished in the same way in classical contact mechanics.

Under a stick condition, no relative motion exists between the bulk and interface systems, so $\tensor{\sigma}_{\text{bulk}} = \tensor{\sigma}_{\text{interface}}$.
Therefore, the stress tensor under a stick condition can be calculated as
\begin{align}
  \tensor{\sigma} = \tensor{\sigma}_{\text{bulk}}\,.
\end{align}
We note that this corresponds to the standard way how a phase-field model of fracture treats a closed crack under compression.
However, existing phase-field models do not allow slip motion along the interface.
This limitation is the major motivation of this work and tackled in the following.

When the interface is under a slip condition, the stress tensor in the interface system is non-zero and different from the bulk stress.
We thus have to evaluate the interface stress such that it incorporates the frictional contact behavior and the no-penetration constraint simultaneously.
For this purpose, we decompose the interface stress tensor into a friction part, $\tensor{\sigma}_{\text{friction}}$, and a no-penetration part, $ \tensor{\sigma}_{\text{no-penetration}}$, as
\begin{align}
  \tensor{\sigma}_{\text{interface}} = \tensor{\sigma}_{\text{friction}} + \tensor{\sigma}_{\text{no-penetration}}\,.
  \label{eq:stress-interface-decomposition}
\end{align}
Of these two, the no-penetration part is determined to make the deformation along the interface and bulk systems compatible in all directions except the slip direction.
This compatibility can be attained by defining the no-penetration part as follows:
\begin{align}
  \tensor{\sigma}_{\text{no-penetration}} = \tensor{\sigma}_{\text{bulk}} - \tau_{\text{bulk}}(\tensor{n}\dyad\tensor{m} + \tensor{m}\dyad\tensor{n})\,,
  \label{eq:stress-no-penetration}
\end{align}
where $\tau_{\text{bulk}} := \tensor{\sigma}_{\text{bulk}}:(\tensor{n}\dyad\tensor{m})$.
In words, the no-penetration part is computed by fully degrading the slip direction components of the bulk stress tensor.
Then, the friction part is calculated to substitute the degraded components, according to a prescribed contact constitutive law.
Under a slip condition, $|\tau| = \mu p_{\cn}$ because the yield function~\eqref{eq:slip-yield-function} is zero, and here $p_{\cn}$ must be equal to $p_{\cn,\text{bulk}} := -\tensor{\sigma}_{\text{bulk}}:(\tensor{n}\dyad\tensor{n})$ because the bulk stress is not degraded in the contact normal direction.
Also, it is likely that the signs of $\tau$ and $\tau_{\text{bulk}}$ are identical.
Therefore, we can express the friction part of the stress tensor as
\begin{align}
  \tensor{\sigma}_{\text{friction}} = \mu p_{\cn,\text{bulk}}\sign(\tau_{\text{bulk}})(\tensor{n}\dyad\tensor{m} + \tensor{m}\dyad\tensor{n})\,.
  \label{eq:stress-friction}
\end{align}
Inserting Eqs.~\eqref{eq:stress-friction} and~\eqref{eq:stress-no-penetration} into Eq.~\eqref{eq:stress-interface-decomposition},
we obtain the interface stress tensor as
\begin{align}
  \tensor{\sigma}_{\text{interface}} = \tensor{\sigma}_{\text{bulk}} + (\mu p_{\cn,\text{bulk}}\sign(\tau_{\text{bulk}}) - \tau_{\text{bulk}})(\tensor{n}\dyad\tensor{m} + \tensor{m}\dyad\tensor{n})\,.
\end{align}
One can see that this stress tensor is obtained by replacing the slip-relevant part of the bulk stress tensor with $\tensor{\sigma}_{\text{friction}}$.
Substituting the above equation into Eq.~\eqref{eq:stress-interpolation} gives the overall stress tensor under a slip condition as
\begin{align}
  \tensor{\sigma} = \tensor{\sigma}_{\text{bulk}} + [1 - g(d)](\mu p_{\cn,\text{bulk}}\sign(\tau_{\text{bulk}}) - \tau_{\text{bulk}})(\tensor{n}\dyad\tensor{m} + \tensor{m}\dyad\tensor{n})\,.
\end{align}
Note that for all contact conditions, we get $\tensor{\sigma} = \tensor{\sigma}_{\text{interface}}$ when $d=1$ and $\tensor{\sigma} = \tensor{\sigma}_{\text{bulk}}$ when $d=0$.

To summarize, we have proposed an approach that calculates the stress tensor in the interface system according to the contact condition of the interface.
In this approach, the interface stress is null under a no-contact condition and equal to the bulk stress under a stick contact condition.
The interface stress under a slip condition is calculated as a combination of the friction part and the no-penetration part so that the contact constitutive behavior and the no-penetration constraints are incorporated into the phase-field formulation.

\smallskip
\begin{remark}
The foregoing expressions for the interface stress tensor can be re-interpreted based on the decomposition of $\tensor{\sigma}_{\text{interface}} = \tensor{\sigma}_{\text{friction}} + \tensor{\sigma}_{\text{no-penetration}}$.
The friction part, $\tensor{\sigma}_{\text{friction}}$, is zero under a no-contact condition, compatible with the bulk stress under a stick condition, and calculated from a friction constitutive law under a slip condition.
The no-penetration part, $\tensor{\sigma}_{\text{no-penetration}}$, is zero under a no-contact condition while it is compatible with the bulk stress under stick and slip contact conditions.
\end{remark}

\section{Discretization and algorithms}
\label{sec:discretization}
This section describes discretization methods and algorithms for numerical solution of the proposed formulation using the standard finite element method.

\subsection{Finite element discretization}
The proposed phase-field formulation can be readily solved by the standard finite element method.
For finite element discretization, we first define trial solution spaces for the displacement field and the phase field as
\begin{align}
  \mathcal{S}_{u} &:= \{\tensor{u} \;\vert\; \tensor{u} \in H^{1}, \; \tensor{u}=\hat{\tensor{u}} \;\;\text{on} \;\; {\pd_{u}\Omega} \}, \\
  \mathcal{S}_{d} &:= \{d \;\vert\; d \in H^{1}\},
\end{align}
where $H^{1}$ denotes a Sobolev space of order one.
Weighting function spaces for the two fields are accordingly defined as
\begin{align}
  \mathcal{V}_{u} &:= \{\tensor{\eta} \;\vert\; \tensor{\eta} \in H^{1}, \; \tensor{\eta}=\tensor{0} \;\;\text{on} \;\; {\pd_{u}\Omega} \}, \\
  \mathcal{V}_{d} &:= \{\phi \;\vert\; \phi \in H^{1} \}\,.
\end{align}
Applying the standard weighted residual procedure, we obtain the following two variational equations:
\begin{align}
  &- \int_{\Omega} \symgrad\tensor{\eta}:\tensor{\sigma}\,\od V
  + \int_{\pd_{t}{\Omega}} \tensor{\eta}\cdot\hat{\tensor{t}}\,\od A = 0\,, \label{eq:var-mom} \\
  &\int_{\Omega} \phi g'(d)\mathcal{H}\, \od V
  + \int_{\Omega} G_{c}\left(\phi\frac{d}{L}\ + L\grad{\phi}\cdot\grad{d}\right) \od V = 0 \,. \label{eq:var-phasefield}
\end{align}
Here, Eq.~\eqref{eq:var-mom} is the linear momentum balance equation
and Eq.~\eqref{eq:var-phasefield} is the phase-field equation.
Both of them are solved in each load step if the interface system is subjected to growth during the course of loading.
However, for a stationary interface problem,
the phase-field equation~\eqref{eq:var-phasefield} only needs to be solved once in the initialization stage of the problem.
The Galerkin and matrix forms of these equations can be developed in a standard manner, so they are omitted for brevity.

As we are considering linear elasticity in this work, the momentum balance equation~\eqref{eq:var-mom} is linear if the contact condition inside the domain is fixed.
This is because the phase-field method has formulated a frictional crack problem as a continuum problem with stiffness that may spatially vary according to the phase field.
However, when the contact condition of a point is subject to change after loading, the problem is incrementally nonlinear.
Therefore, it is necessary to apply a nonlinear solution method for the momentum equation.
Note that the solution method will converge quickly whenever the contact condition remains unchanged from an initial guess.
The phase-field equation~\eqref{eq:var-phasefield} is always linear so it can be solved easily.

In this work, we use Newton's method to solve the discretized momentum balance equation.
During Newton iterations, the increment of the nodal displacement vector, denoted by $\Del{\tensor{U}}$, can be obtained by solving
\begin{align}
  \tensor{\mathcal{R}} = -\tensor{\mathcal{J}}\Del{\tensor{U}}\,,
\end{align}
where $\tensor{\mathcal{R}}$ is the residual vector, of which element-wise contribution can be calculated as
\begin{align}
  [\tensor{\mathcal{R}}]_{e}^{i} :=
  - \int_{\Omega_e} \symgrad\tensor{\eta}^{i}:\tensor{\sigma}\,\od V
  + \int_{\pd_{t}{\Omega}_{e}} \tensor{\eta}^{i}\cdot\hat{\tensor{t}}\,\od A\,,
  \label{eq:residual}
\end{align}
and $\tensor{\mathcal{J}}$ is the Jacobian matrix, of which element-wise contribution can be calculated as
\begin{align}
  [\tensor{\mathcal{J}}]_{e}^{i,j} :=
  - \int_{\Omega_e} \symgrad\tensor{\eta}^{i}:\mathbb{C}:\symgrad\tensor{\eta}^{j}\,\od V\,,
  \label{eq:jacobian}
\end{align}
with $e$ denoting an element index and $i,j$ denoting shape function indices.
$\mathbb{C}$ is the stress--strain tangent that is the same as $\mathbb{C}_{\text{bulk}}$ in the bulk region
but may be different from it otherwise.

\smallskip
\begin{remark}
\revised{
As is well known, phase-field modeling requires very fine discretization around where a discontinuity is approximated diffusely. This requirement inevitably makes the computation cost for running the phase-field method more expensive than that for discrete methods such as AES and XFEM.
However, because the phase-field method can be implemented far more easily than these discrete methods,
the total cost for implementing and running the phase-field method would be attractively low.
Note that a local and/or adaptive mesh refinement can reduce the running cost significantly.
Also, numerical solutions to phase-field formulations are known to be insensitive to mesh alignment, see Mandal \etal~\cite{Mandal2019} for example.
Therefore, both structured and unstructured meshes can be used as long as their element sizes are small enough around the interface region.
The question of how small is enough for element sizes will be addressed through numerical examples in Section~\ref{sec:examples}.
}
\end{remark}

\subsection{Update of stress and tangent tensors}
To evaluate Eqs.~\eqref{eq:residual} and~\eqref{eq:jacobian} during finite element assembly,
we need to calculate the stress tensor, $\tensor{\sigma}$, and the stress--strain tangent tensor, $\mathbb{C}$, at every quadrature point.
Consider a typical Newton update step for which the strain tensor, $\tensor{\eps}$, and the phase-field variable, $d$, are given
at the quadrature point of interest.
It is also assumed that the unit vector in the interface normal direction, $\tensor{n}$ and the unit vector along the slip direction, $\tensor{m}$, are also known at this point.

The stress tensor at a quadrature point can be updated as described in Algorithm~\ref{algo:stress-update}.
The algorithm first checks whether the current quadrature point belongs to a bulk region where $d=0$.
If not, the algorithm identifies the contact condition at the quadrature point and then updates the stress tensor following the approach proposed in Section~\ref{subsec:stress-calculation}.
To distinguish between stick and slip conditions, the algorithm uses a standard predictor--corrector approach employing the bulk stress as a trial stress.
Therefore, the yield function, $f$, is calculated using $\tau$ and $p_{\cn}$ of the bulk stress.
\begin{algorithm}[htbp]
  \caption{Stress update procedure for phase-field modeling of frictional cracks}
  \begin{algorithmic}[1]
  \Require $\tensor{\eps}$, $d$, $\tensor{n}$, and $\tensor{m}$ as well as material parameters at a quadrature point.
  \If {$d = 0$}
    \State Bulk region. Return $\tensor{\sigma}=\mathbb{C}_{\text{bulk}}:\tensor{\eps}$.
  \EndIf
  \State Calculate the interface normal strain $\eps_{N}=\tensor{\eps}:(\tensor{n}\dyad\tensor{n})$ and the bulk stress $\tensor{\sigma}_{\text{bulk}}=\mathbb{C}_{\text{bulk}}:\tensor{\eps}$.
  \If {$\eps_{\cn} > 0$}
    \State No-contact condition. Return $\tensor{\sigma} = g(d)\tensor{\sigma}_{\text{bulk}}$.
  \EndIf
  \State Calculate the yield function by using the bulk stress $f = |\tau_{\text{bulk}}| - \mu p_{\cn,\text{bulk}}$,
  where $\tau_{\text{bulk}} = \tensor{\sigma}_{\text{bulk}}:(\tensor{n}\dyad\tensor{m})$
  and $p_{\cn,\text{bulk}} = -\tensor{\sigma}_{\text{bulk}}:(\tensor{n}\dyad\tensor{n})$.
  \If {$f < 0$}
    \State Stick condition. Return $\tensor{\sigma}=\tensor{\sigma}_{\text{bulk}}$.
  \Else
    \State Slip condition. Return $\tensor{\sigma} = \tensor{\sigma}_{\text{bulk}} + [1 - g(d)](\mu p_{\cn,\text{bulk}}\sign(\tau_{\text{bulk}}) - \tau_{\text{bulk}})(\tensor{n}\dyad\tensor{m} + \tensor{m}\dyad\tensor{n})$.
  \EndIf
  \Ensure $\tensor{\sigma}$ at the quadrature point.
  \end{algorithmic}
  \label{algo:stress-update}
\end{algorithm}

The stress--strain tangent tensor should also be evaluated to assemble the Jacobian matrix, see Eq.~\eqref{eq:jacobian}.
This calculation is trivial for a bulk region.
For an interface region where $0<d\leq1$, it is given by
\begin{align}
  \mathbb{C} =
  \left\{\begin{array}{ll}
  g(d)\mathbb{C}_{\text{bulk}} & \text{for a no-contact condition,}\\
  \mathbb{C}_{\text{bulk}} & \text{for a stick condition,}\\
  \mathbb{C}_{\text{bulk}} + [1-g(d)](\mathbb{C}_{f} - \mathbb{C}_{\tau}) & \text{for a slip condition.}
  \end{array}\right.
\end{align}
Here, for a slip condition, $\mathbb{C}_{f}$ is defined as
\begin{align}
  \mathbb{C}_{f} &:= -\sign(\tau_{\text{bulk}})\mu[\lambda(\tensor{n}\dyad\tensor{m} + \tensor{m}\dyad\tensor{n})\dyad\tensor{1}
  + 2G(\tensor{n}\dyad\tensor{m} + \tensor{m}\dyad\tensor{n})\dyad(\tensor{n}\dyad\tensor{n})]\,,
\end{align}
and $\mathbb{C}_{\tau}$ is defined as
\begin{align}
  \mathbb{C}_{\tau} &:= G[(\tensor{n}\dyad\tensor{m} + \tensor{m}\dyad\tensor{n})\dyad(\tensor{n}\dyad\tensor{m} + \tensor{m}\dyad\tensor{n})]\,.
\end{align}
Note that all these tensors can be calculated in a straightforward manner
as long as the vectors $\tensor{n}$ and $\tensor{m}$ are given.

\subsection{Calculation of unit normal and slip vectors}
As described above, the unit vector in the interface normal direction, $\tensor{n}$,
and the unit vector in the slip direction, $\tensor{m}$,
are crucial for the proposed phase-field formulation.
Usually, in phase-field modeling, the interface normal vector is approximated as $\tensor{n}\approx\grad{d}/\|\grad{d}\|$.
The accuracy of this approximation, however, seems to be insufficient for our purpose, for two main reasons:
(1) because a crack tip has been approximated bluntly, $\grad{d}/\|\grad{d}\|$ calculated around the crack tip region
is indeed nearly orthogonal to the desired $\tensor{n}$, and
(2) unless $d$ is very close to 1, $\grad{d}/\|\grad{d}\|$ may not be close enough to the desired $\tensor{n}$.

Therefore, for a more accurate calculation of $\tensor{n}$ and $\tensor{m}$, we devise a new algorithm that first identifies a lower-dimensional crack path from phase-field values and then estimates the unit vectors from the identified crack path.
Hereafter, we will restrict our attention to 2D problems because identifying a crack path in a 3D phase field is an outstanding challenge.
In a 2D phase-field problem, a crack path would be a 1D line that connects points where $d=1$.
Drawing on this idea, we construct such a line through the procedure described in Algorithm~\ref{algo:normal-slip-vectors}.
In essence, this algorithm finds nodes where $d\approx 1$ that are distant at least the length parameter $L$
and connect them in a piecewise linear manner to approximate the crack path.
From this piecewise linearly approximated crack path, we can compute $\tensor{n}$ and $\tensor{m}$ and then assign them to nearby quadrature points.
\begin{algorithm}[htbp]
  \caption{Calculation of unit normal and slip vectors for a phase-field approximated interface}
  \begin{algorithmic}[1]
  \Require Coordinates and phase-field values of nodes.
  \State Find nodes where phase-field values are greater than a threshold (\eg~0.98) and store them to a set $\Gamma_{\text{tmp}}$.
  \State Search the node $\mathcal{N}_{1}$ where the phase field value is greatest among all nodes in $\Gamma_{\text{tmp}}$.
  \State Find all nodes whose distances from $\mathcal{N}_{1}$ are within the phase-field length parameter, $L$; then remove them from $\Gamma_{\text{tmp}}$ and move $\mathcal{N}_{1}$ to a new set $\Gamma_{\text{final}}$.
  \State Repeat the above step for other nodes in the order of decreasing phase-field value, say $\mathcal{N}_{2}, \mathcal{N}_{3}, \cdots$, until $\Gamma_{\text{tmp}}$ becomes empty.
  \State Sort nodes in $\Gamma_{\text{final}}$ according to their coordinates in one direction (e.g. the $x$ direction).
  \State Connected the sorted nodes in $\Gamma_{\text{final}}$ as a piecewise linear line. The piecewise linear line is then considered a crack path.
  \State For each segment in the piecewise linear crack path, calculate $\tensor{n}$ and $\tensor{m}$.
  \State For quadrature points where $d>0$, assign $\tensor{n}$ and $\tensor{m}$ from the nearest segment in the piecewise linear crack path.
  \Ensure $\tensor{n}$ and $\tensor{m}$ at quadrature points where $d>0$.
  \end{algorithmic}
  \label{algo:normal-slip-vectors}
\end{algorithm}

The proposed algorithm for calculating $\tensor{n}$ and $\tensor{m}$ is simple and appears to be sufficiently accurate for phase-field modeling of frictional interfaces.
We note, however, that any other algorithm can be used for the same purpose as long as it gives reliable results for the unit vectors.
For example, Ziaei-Rad \etal~\cite{ZiaeiRad2016} have proposed a variational method for identifying a crack path in phase-field modeling.
The use of such an advanced method will likely improve the accuracy of the overall numerical solution, although it requires significantly more effort for implementation.
Furthermore, due to lack of a good algorithm for estimating $\tensor{n}$ and $\tensor{m}$ in 3D, applications of the proposed phase-field formulation will be limited to 2D problems in the next section.
Overcoming this limitation for 3D problems will be a future research topic.

\section{Numerical examples}
\label{sec:examples}
This section has two objectives: (1) to verify the proposed phase-field formulation for frictional interfaces,
and (2) to demonstrate the capability of the proposed method for modeling crack growth with frictional contact.
For the first objective, we adopt three numerical examples of frictional interfaces that have been simulated by discrete methods in the literature.
Yet the interfaces in these benchmark examples are stationary (\ie~not allowed to advance).
Therefore, for the second objective, we introduce a fourth example whereby a preexisting crack propagates according to the phase-field equation~\eqref{eq:var-phasefield}.

Results in this section have been obtained using the \texttt{deal.II} finite element library~\cite{Bangerth2007,dealII90}.
Bilinear quadrilateral elements have been used for all numerical examples.
Plane strain conditions are assumed throughout.

\subsection{Square domain with an internal crack}
\begin{figure}[htbp]
  \centering
  \includegraphics[width=0.55\textwidth]{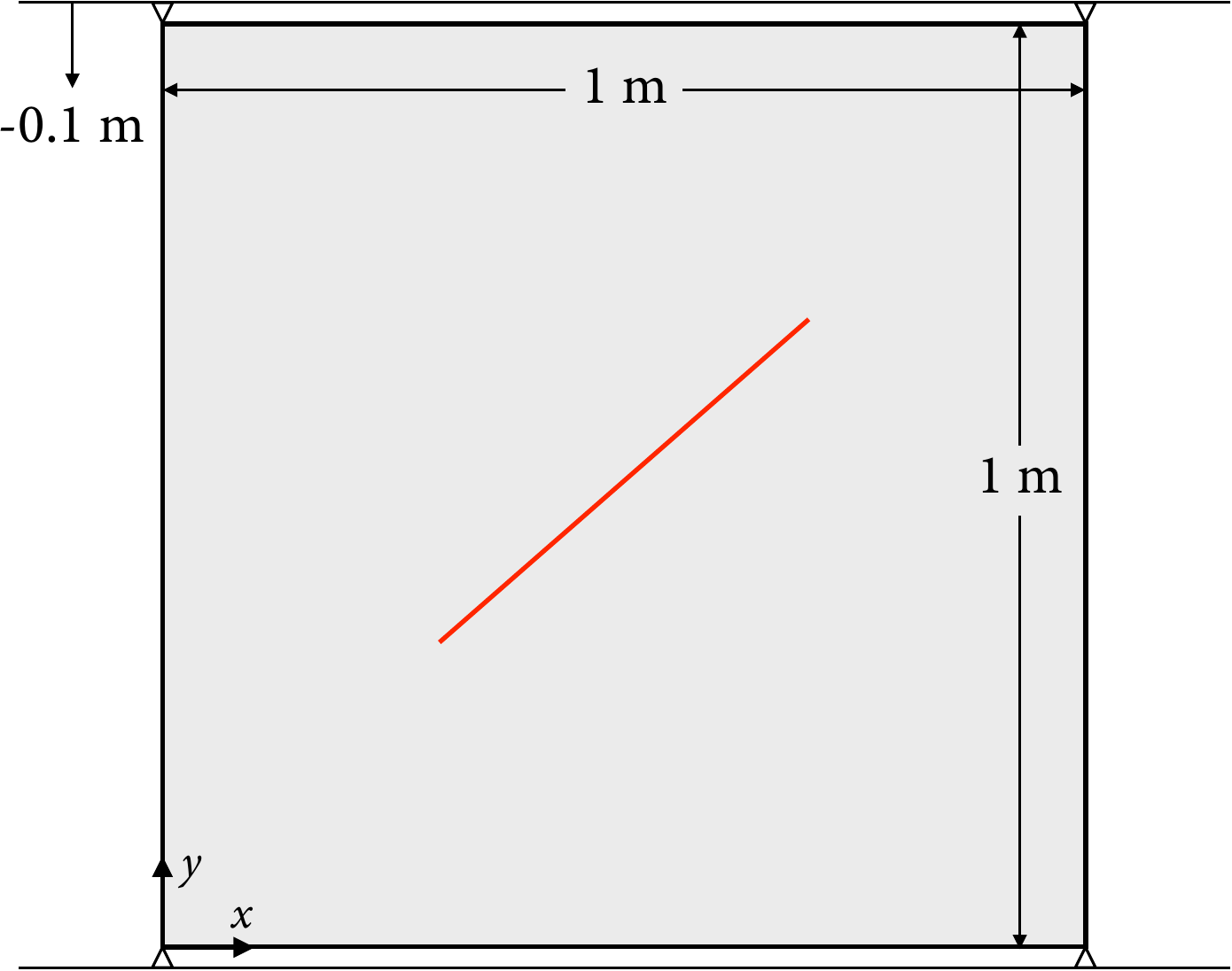}
  \caption{Setup of the internal crack problem.}
  \label{fig:internal-crack-setup}
\end{figure}

Our first example is compression of an internally cracked domain depicted in Fig.~\ref{fig:internal-crack-setup}.
This problem was initially presented in Dolbow \etal~\cite{Dolbow2001} and then used by other studies including Liu and Borja~\cite{Liu2008} and Annavarapu \etal~\cite{Annavarapu2014}.
Here we also consider this problem for verification of the phase-field formulation.
The domain is a 1 m wide square and possesses a crack whose tips are located at coordinates $(0.3, 0.33)$ m and $(0.7, 0.68)$ m.
Note that these tip locations are adopted from Annavarapu \etal~\cite{Annavarapu2014} and slightly different from those in Liu and Borja~\cite{Liu2008}.
We assign the elasticity parameters of the material as $E=10000$ MPa and $\nu=0.3$,
and the friction coefficient of the crack as $\mu=0.1$.
The top boundary of the domain is subjected to a prescribed displacement of $-0.1$ m (downward).

To investigate the convergence of numerical solution with the length parameter and the element size,
we consider three values of $L$, namely $0.008$ m, $0.004$ m, and $0.002$ m and three values of $L/h$, namely 4, 8, and 16.
The values of $L/h$ are selected based on results in Borden \etal~\cite{Borden2012} that indicate $L/h \geq 4$ gives reasonably accurate solutions when the same type of surface density functional is employed.
Because such fine discretization is necessary only for the interface region and its nearby, we locally refine elements around a node where the phase-field variable is greater than a threshold, until their size reaches a prescribed $L/h$ value.
To determine the threshold value, we recall that the spatial variation of the phase field in the chosen surface density functional is given by $d=\exp(-|x|/L)$, with $x$ denoting the distance from the point where $d=1$~\cite{Miehe2010}.
As $d=\exp(-1)\approx0.378$ when $x=L$, we set the threshold as 0.1 to make the locally refined region sufficiently wide.
The same mesh refinement scheme will be used throughout this section.

\begin{figure}[htbp]
  \centering
  \includegraphics[width=0.95\textwidth]{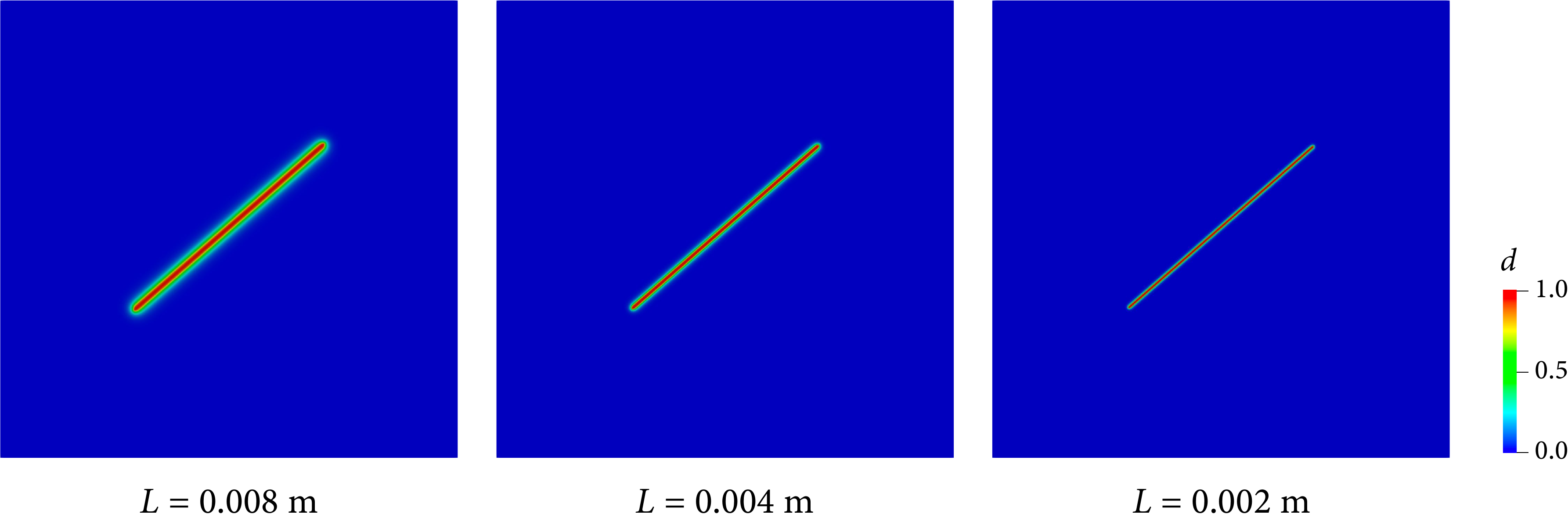}
  \caption{Phase-field distributions initialized to diffusely approximate the internal crack depicted in Fig.~\ref{fig:internal-crack-setup}.}
  \label{fig:internal-crack-pf}
\end{figure}

For initialization of the phase-field variable, we adopt a standard way in phase-field fracture modeling that prescribes $\mathcal{H}$ at quadrature points around a preexisting crack (see Appendix A of Borden \etal~\cite{Borden2012} for example).
With $\mathcal{H}$ values prescribed to make $d=1$ at the initial crack, we solve the phase-field equation~\eqref{eq:var-phasefield} once to obtain a phase-field distribution that will be used throughout the problem.
Figure~\ref{fig:internal-crack-pf} shows phase-field distributions in the $L=0.008$ m, $0.004$ m, and $0.002$ m cases when $L/h=8$.
It is clear that the diffuse approximation zone becomes narrower as $L$ decreases.
After initializing the phase field in this way, we simulate the problem through 10 load steps with a uniform displacement increment of -0.01 m on the top boundary.

\begin{figure}[htbp]
  \centering
  \subfloat[$x$-displacement\vspace{1em}]{\includegraphics[width=0.95\textwidth]{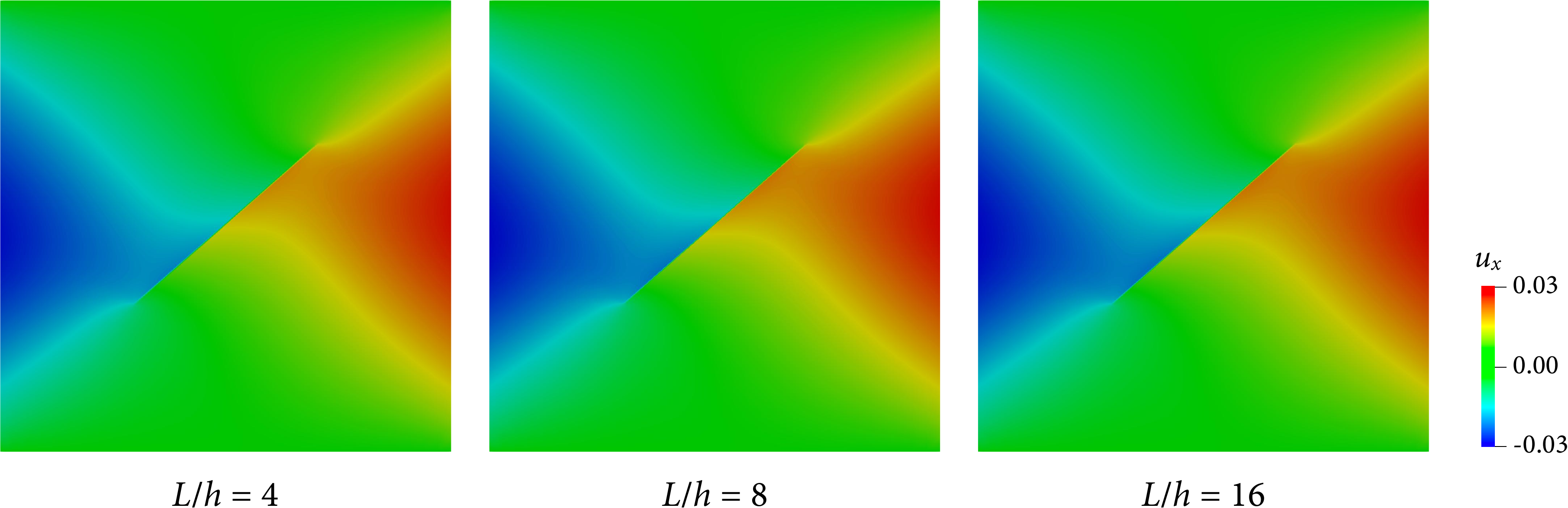}}\\ \vspace{1em}
  \subfloat[$y$-displacement\vspace{1em}]{\includegraphics[width=0.95\textwidth]{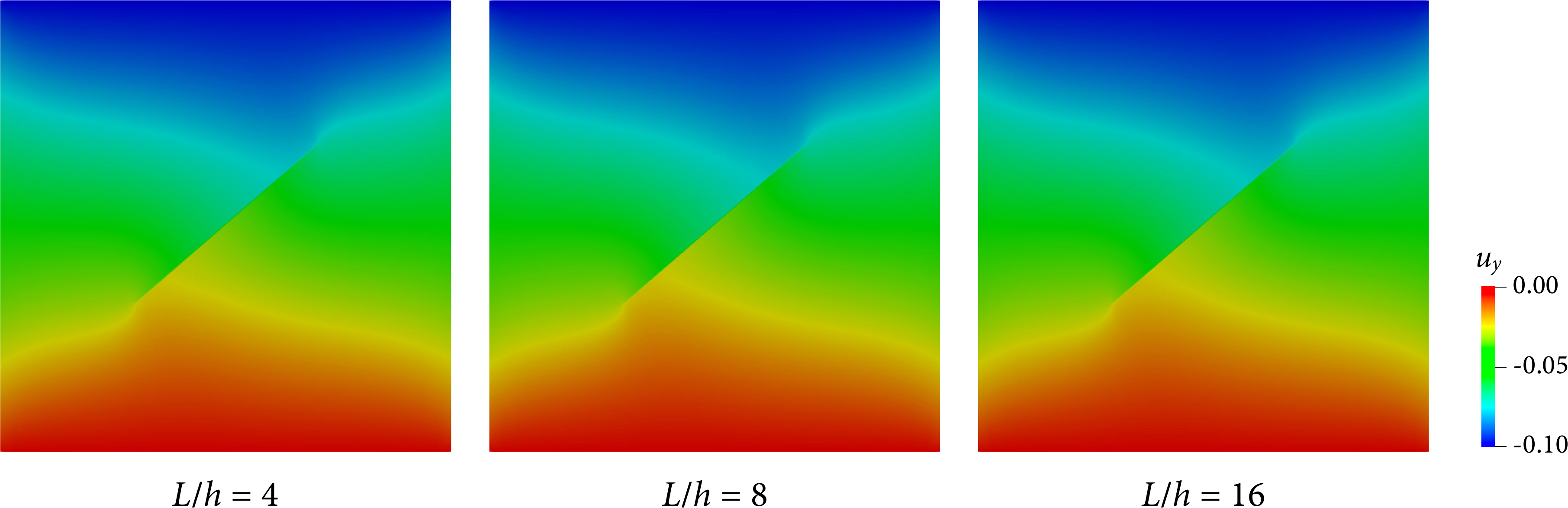}}
  \caption{Displacement fields in mesh refinement tests with $L=0.002$ m. Color bar in meters.}
  \label{fig:internal-crack-h-refinement}
\end{figure}

\begin{figure}[htbp]
  \centering
  \subfloat[Normal pressures]{\includegraphics[height=0.265\textheight]{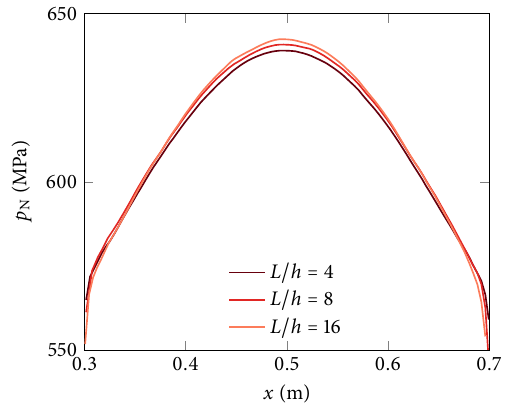}}$\,\,$
  \subfloat[Tangential stresses]{\includegraphics[height=0.265\textheight]{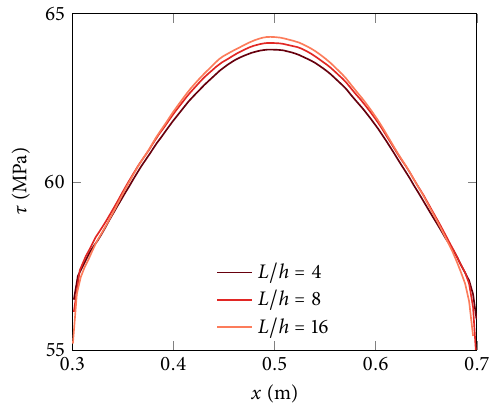}}
  \caption{Contact normal pressures and tangential stresses in mesh refinement tests with $L=0.002$ m.}
  \label{fig:internal-crack-h-refinement-stresses}
\end{figure}

We begin by checking to see whether numerical solutions converge with mesh refinement, \revised{in terms of displacement fields in the domain and contact normal pressures and tangential stresses along the crack.
The contact pressures and stresses are obtained at the nodes in $\Gamma_{\text{final}}$ (see Algorithm~\ref{algo:normal-slip-vectors}) after nodal projection.
Figures~\ref{fig:internal-crack-h-refinement} and~\ref{fig:internal-crack-h-refinement-stresses} present these results from a mesh refinement study conducted for $L=0.002$ m.
We observe that both the displacement and stress solutions converge as the element size decreases.
Note also that the contact pressures and stresses satisfy the prescribed contact law ($\tau = 0.1 p_{\cn}$ here) in every mesh.}
Although not presented, results of mesh refinement studies conducted with other length parameters showed more or less the same patterns.
Therefore, we have found that the numerical model converges with mesh refinement, and that refinement levels of $L/h \geq 4$ provide reasonable accuracy.

\begin{figure}[htbp]
  \centering
  \subfloat[$x$-displacement\vspace{1em}]{\includegraphics[width=0.95\textwidth]{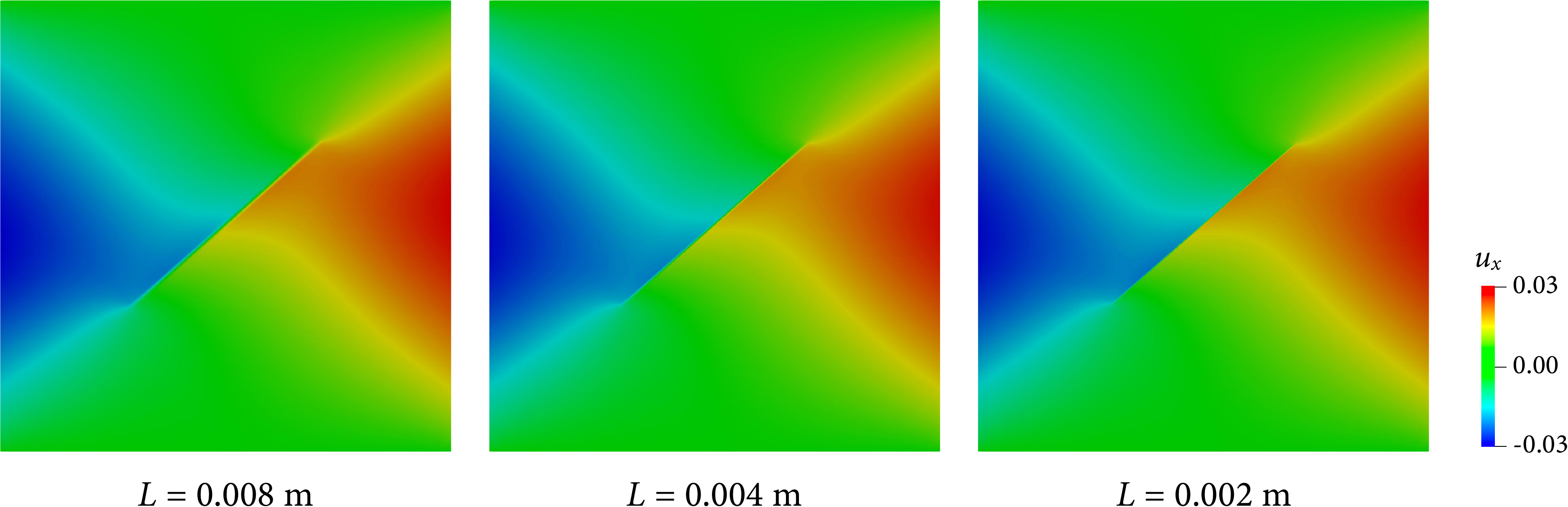}}\\ \vspace{1em}
  \subfloat[$y$-displacement\vspace{1em}]{\includegraphics[width=0.95\textwidth]{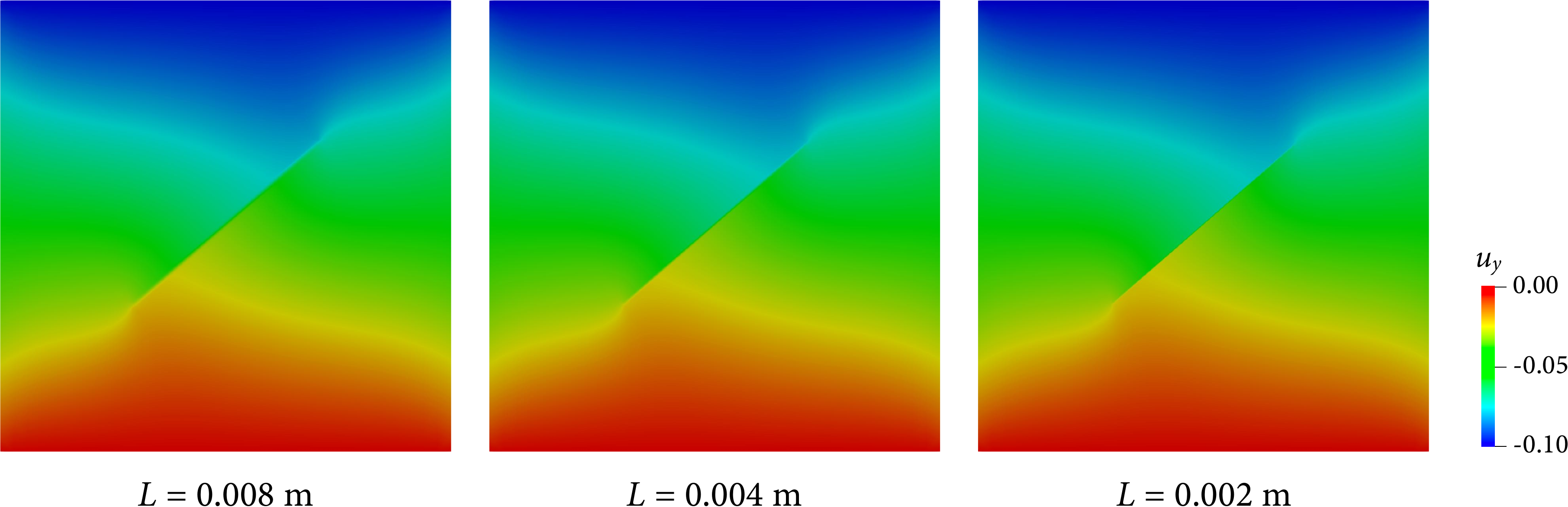}}
  \caption{Displacement fields in length parameter sensitivity tests with $L/h=8$. Color bar in meters.}
  \label{fig:internal-crack-L-refinement}
\end{figure}

\begin{figure}[htbp]
  \centering
  \subfloat[Normal pressures]{\includegraphics[height=0.265\textheight]{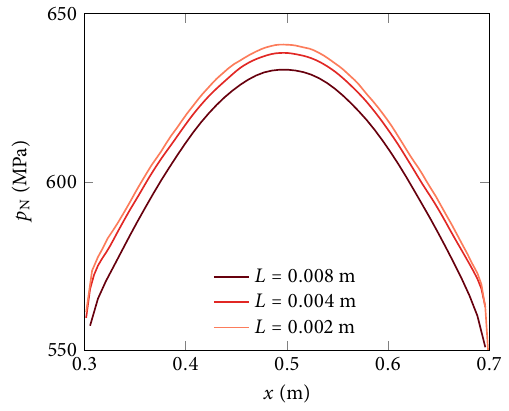}}$\,\,$
  \subfloat[Tangential stresses]{\includegraphics[height=0.265\textheight]{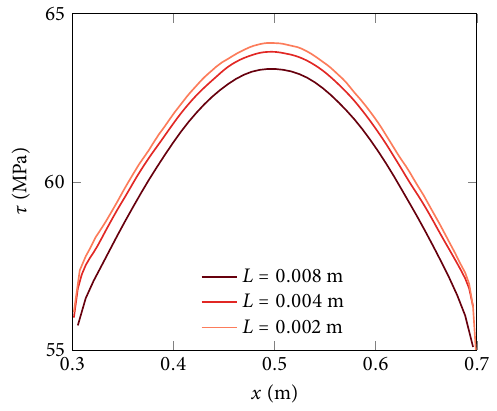}}
  \caption{Contact normal pressures and tangential stresses in length parameter sensitivity tests with $L/h=8$.}
  \label{fig:internal-crack-L-refinement-stresses}
\end{figure}

Next, we examine how numerical solutions are sensitive to the length parameter for phase-field approximation.
\revised{Figures~\ref{fig:internal-crack-L-refinement} and~\ref{fig:internal-crack-L-refinement-stresses} show displacement fields and contact pressures and stresses, respectively, obtained using the three length parameters with meshes of $L/h=8$.}
One can see that the results show little sensitivity to $L$.
Although a smaller $L$ leads to a sharper displacement jump across the crack and thereby larger contact stresses, the $L=0.008$ m case also shows fairly good results.
This observation indicates that a rather diffuse approximation still can provide reasonably good solutions.
Therefore, it can be concluded that the phase-field method does not require a very small $L$ for practical purposes while allowing us to obtain a more accurate solution by reducing $L$.

\begin{figure}[htbp]
  \centering
  \includegraphics[width=0.85\textwidth]{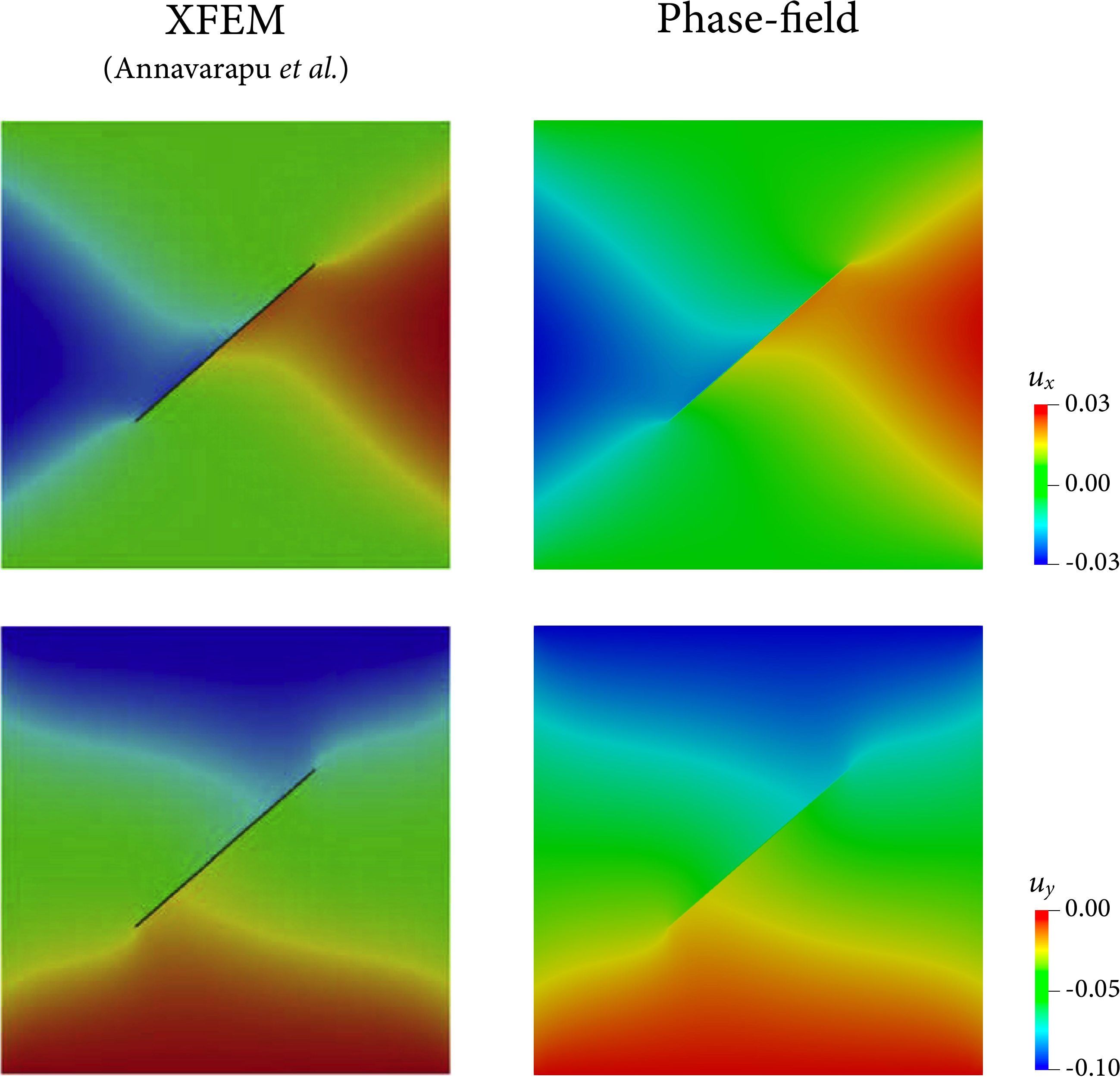}
  \caption{Comparison of phase-field solutions with XFEM solutions in Annavarapu \etal~\cite{Annavarapu2014}. The phase-field results have been obtained with $L=0.002$ m and $L/h=8$. Color bar in meters.}
  \label{fig:internal-crack-comparison}
\end{figure}

Having confirmed that the proposed method gives consistent solutions, we now verify it with results in the literature.
We particularly compare our results from the $L=0.002$ m and $L/h=8$ case with results in Annavarapu \etal~\cite{Annavarapu2014} obtained by the combination of XFEM and a weighted Nitsche method.
Figure~\ref{fig:internal-crack-comparison} shows this comparison.
As can be seen, the phase-field and XFEM results are nearly identical in both qualitative and quantitative aspects.
Therefore, we have verified that the proposed phase-field method can provide numerical solutions comparable to those obtained by advanced discrete methods for frictional cracks.

Lastly, we remark that the phase-field formulation for this problem is linear in all the load steps, requiring only a single Newton update.
This is because the contact condition along the crack is identified as a slip condition from the initial stress-free condition (as zero stress makes $f=0$), and it is indeed a slip condition throughout the problem.
So this example was nothing but a linear elasticity problem with heterogeneous stiffness.
This means that, although the phase-field method requires a quite fine mesh,
its fast convergence can counterbalance the overall computational cost.

\subsection{Square domain with an inclined interface}
\begin{figure}[htbp]
  \centering
  \includegraphics[width=0.55\textwidth]{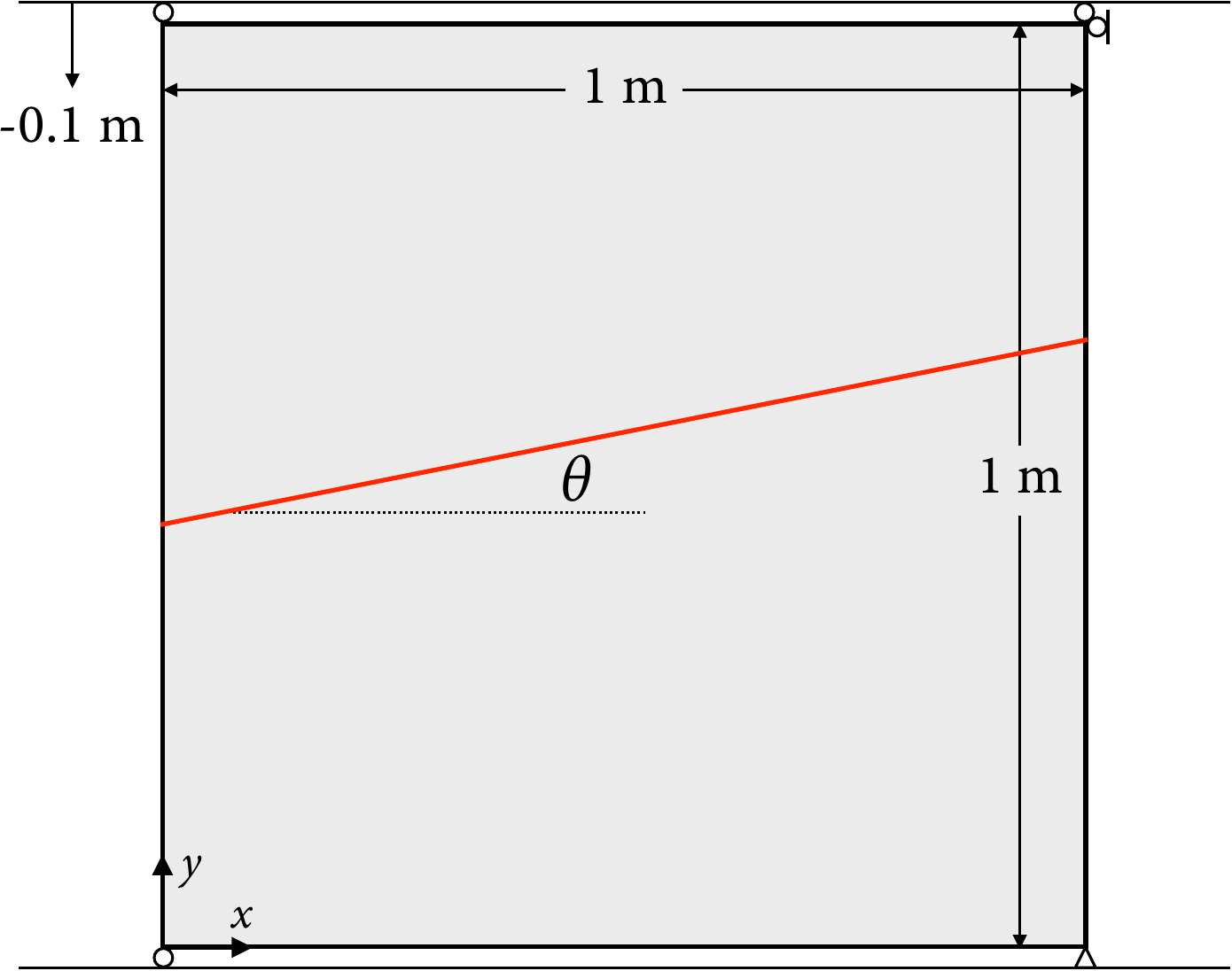}
  \caption{Setup of the inclined interface problem.}
  \label{fig:inclined-interface-setup}
\end{figure}

The purpose of our second example is to investigate the ability of the proposed phase-field method to distinguish between stick and slip conditions.
For this purpose, we adopt the problem of a square domain with an inclined interface, which was also first used in Dolbow \etal~\cite{Dolbow2001} and later revisited by Annavarapu \etal~\cite{Annavarapu2014}, among others.
The setup of this problem is illustrated in Fig.~\ref{fig:inclined-interface-setup}.
Similar to the previous example, a 1 m wide square domain is compressed from the top, but here the discontinuous interface is extended to the side boundaries of the domain.
The interface is inclined from the horizontal with an angle $\theta=\tan^{-1}(0.2)$.
Therefore, when the friction coefficient $\mu$ is smaller than 0.2, the upper block should slip along the interface; otherwise, the upper and lower blocks should be sticked together and behave as a whole.
Accordingly, this problem serves as a good benchmark example for examining the capability for distinguishing between stick and slip behaviors.

Because the domain size remains the same as the previous example, we consider the same three length parameters, $L=0.008$ m, 0.004 m, and 0.002 m.
We discretize the domain and initialize the phase field using the same way in the previous example.
The refinement level is now fixed as $L/h=8$.
Following Annavarapu \etal~\cite{Annavarapu2014}, we consider two cases of friction coefficients, namely $\mu=0.19$ and $\mu=0.21$.
The elasticity parameters are set as $E=1000$ MPa and $\nu=0.3$ for both the upper and lower blocks.
We again use 10 load steps with a uniform displacement increment of $-0.01$ m.

\begin{figure}[htbp]
  \centering
  \includegraphics[width=0.95\textwidth]{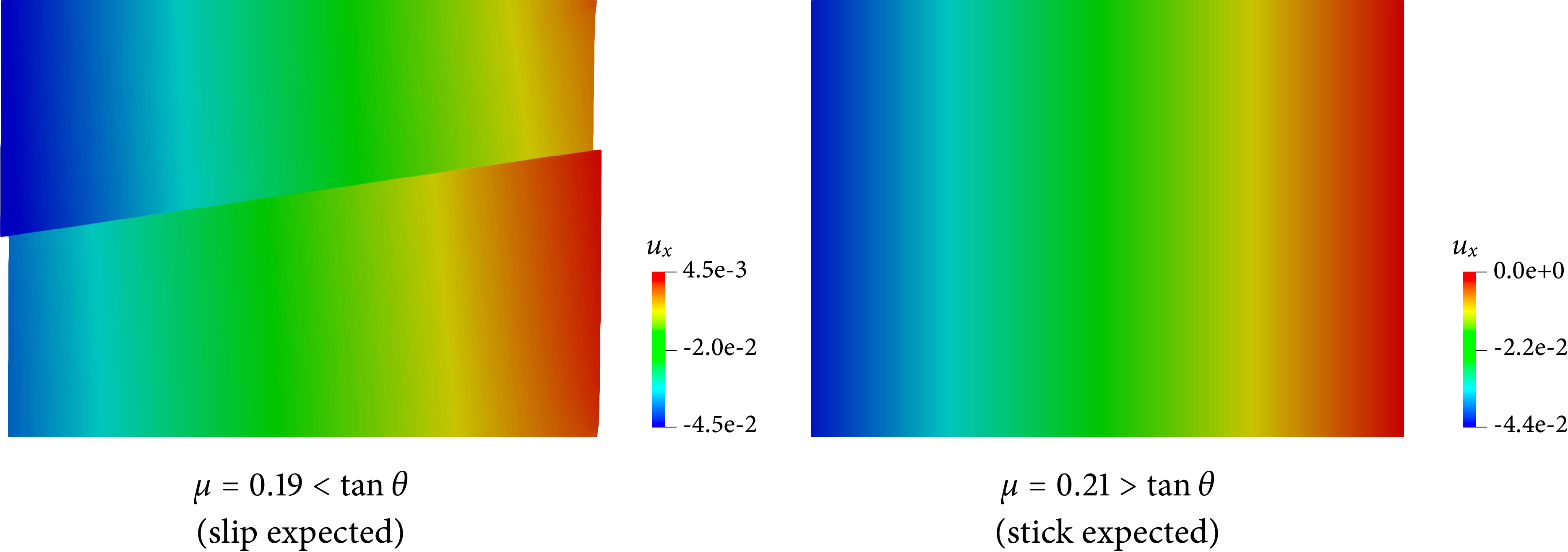}
  \caption{Results of stick/slip distinction tests with $L=0.002$ m. Color bar in meters. Displacement is scaled by a factor of 2.}
  \label{fig:inclined-interface-stick-slip}
\end{figure}

Figure~\ref{fig:inclined-interface-stick-slip} comparatively shows the results of $\mu=0.19$ and $\mu=0.21$ cases in terms of the $x$-displacement field.
We can find that the domain is under a slip condition when $\mu=0.19<\tan\theta$ and under a stick condition when $\mu=0.21>\tan\theta$.
Therefore, we have confirmed that the phase-field method can also distinguish between stick and slip conditions appropriately.

\begin{figure}[htbp]
  \centering
  \includegraphics[width=0.95\textwidth]{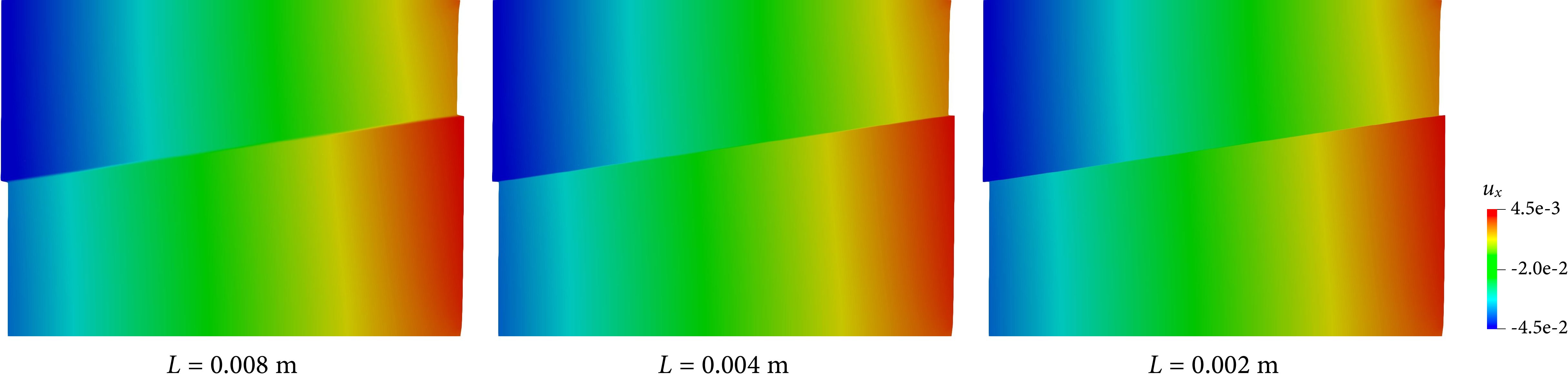}
  \caption{Displacement fields in length parameter sensitivity tests with $L/h=8$. Color bar in meters. Displacement is scaled by a factor of 2.}
  \label{fig:inclined-interface-L-refinement}
\end{figure}

\begin{figure}[htbp]
  \centering
  \subfloat[Normal pressures]{\includegraphics[height=0.265\textheight]{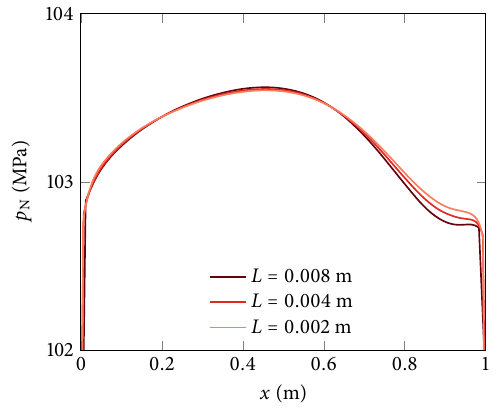}}$\,\,$
  \subfloat[Tangential stresses]{\includegraphics[height=0.265\textheight]{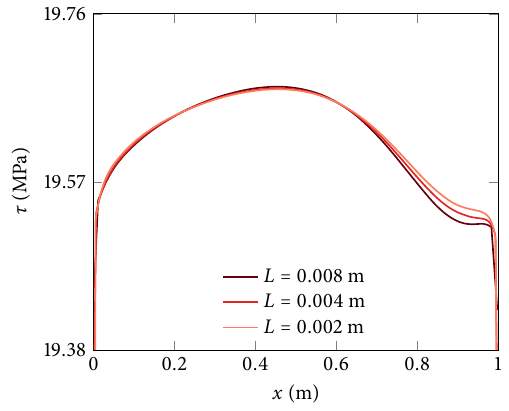}}
  \caption{Contact normal pressures and tangential stresses in length parameter sensitivity tests with $L/h=8$.}
  \label{fig:inclined-interface-L-refinement-stresses}
\end{figure}

Also for this problem, we check the sensitivity to the length parameter by repeating the same problem with $L=0.008$ m, 0.004 m, and 0.002 m.
We have found that stick and slip conditions are correctly distinguished with all the three length parameters.
Because the stick case is a standard linear elasticity problem, we only present the displacement and contact stress results of the slip case ($\mu=0.19$) in Figs.~\ref{fig:inclined-interface-L-refinement} and \ref{fig:inclined-interface-L-refinement-stresses}.
It can again be seen that the numerical solutions show little sensitivity to $L$, at least for the length parameters considered which are reasonably small compared with the domain size.
\revised{Results in Fig.~\ref{fig:inclined-interface-L-refinement-stresses} also confirm that the numerical solutions converge as $L$ decreases.}

\begin{figure}[htbp]
  \centering
  \subfloat[$\mu=0.19$ (slip) -- initial condition: slip]{\includegraphics[width=0.49\textwidth]{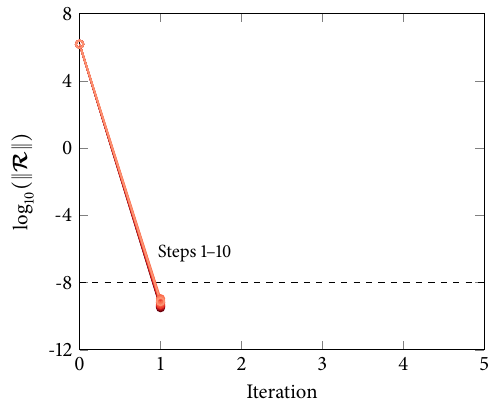}}$\,\,$
  \subfloat[$\mu=0.21$ (stick)  -- initial condition: slip]{\includegraphics[width=0.49\textwidth]{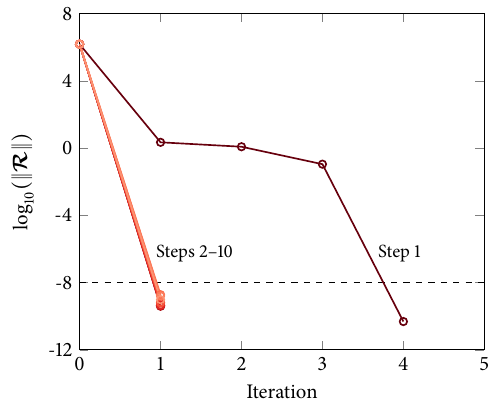}}
  \caption{Residual vector norms during Newton iterations in the $\mu=0.19$ (slip) and $\mu=0.21$ (stick) cases, when the initial contact condition is a slip condition.}
  \label{fig:inclined-interface-newton-slip-guess}
\end{figure}

\begin{figure}[htbp]
  \centering
  \subfloat[$\mu=0.19$ (slip) -- initial condition: stick]{\includegraphics[width=0.49\textwidth]{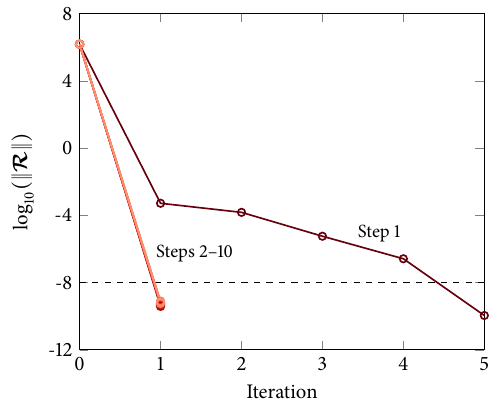}}$\,\,$
  \subfloat[$\mu=0.21$ (stick) -- initial condition: stick]{\includegraphics[width=0.49\textwidth]{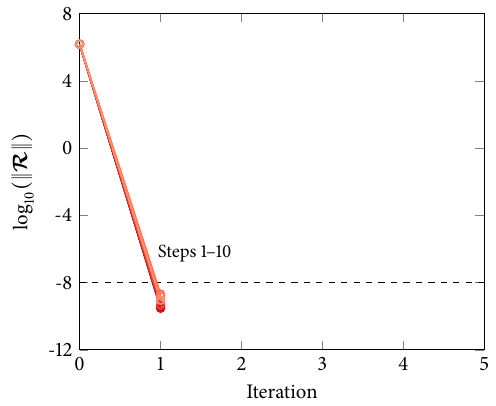}}
  \caption{Residual vector norms during Newton iterations in the $\mu=0.19$ (slip) and $\mu=0.21$ (stick) cases, when the initial contact condition is a stick condition.}
  \label{fig:inclined-interface-newton-stick-guess}
\end{figure}

Figure~\ref{fig:inclined-interface-newton-slip-guess} shows the Newton convergence behaviors of the stick and slip cases when $L=0.002$ m and $L/h=8$.
We can see that except the first load step of the stick case, all load steps converged after a single update, which evinces the linearity of the formulation.
The first step of the stick case required multiple iterations because a slip condition is initially assumed for the stress-free initial condition and it has to be corrected by a Newton iteration.
From the second step, as a stick condition is identified from the last converged step, the problem remained linear.
To confirm this statement, we have also repeated the same problems by changing the initial contact condition to a stick condition.
Then, as shown in Fig.~\ref{fig:inclined-interface-newton-stick-guess}, the first load step of the slip case required two iterations for convergence, and all other load steps in the stick and slip cases converged after a single Newton update.
Therefore, we can conclude that the formulation is linear if the initial guess of the contact condition is correct, and that
an incorrect guess of the contact condition can be rectified during a Newton iteration.

\subsection{Sliding of a block}
\begin{figure}[htbp]
  \centering
  \includegraphics[width=0.8\textwidth]{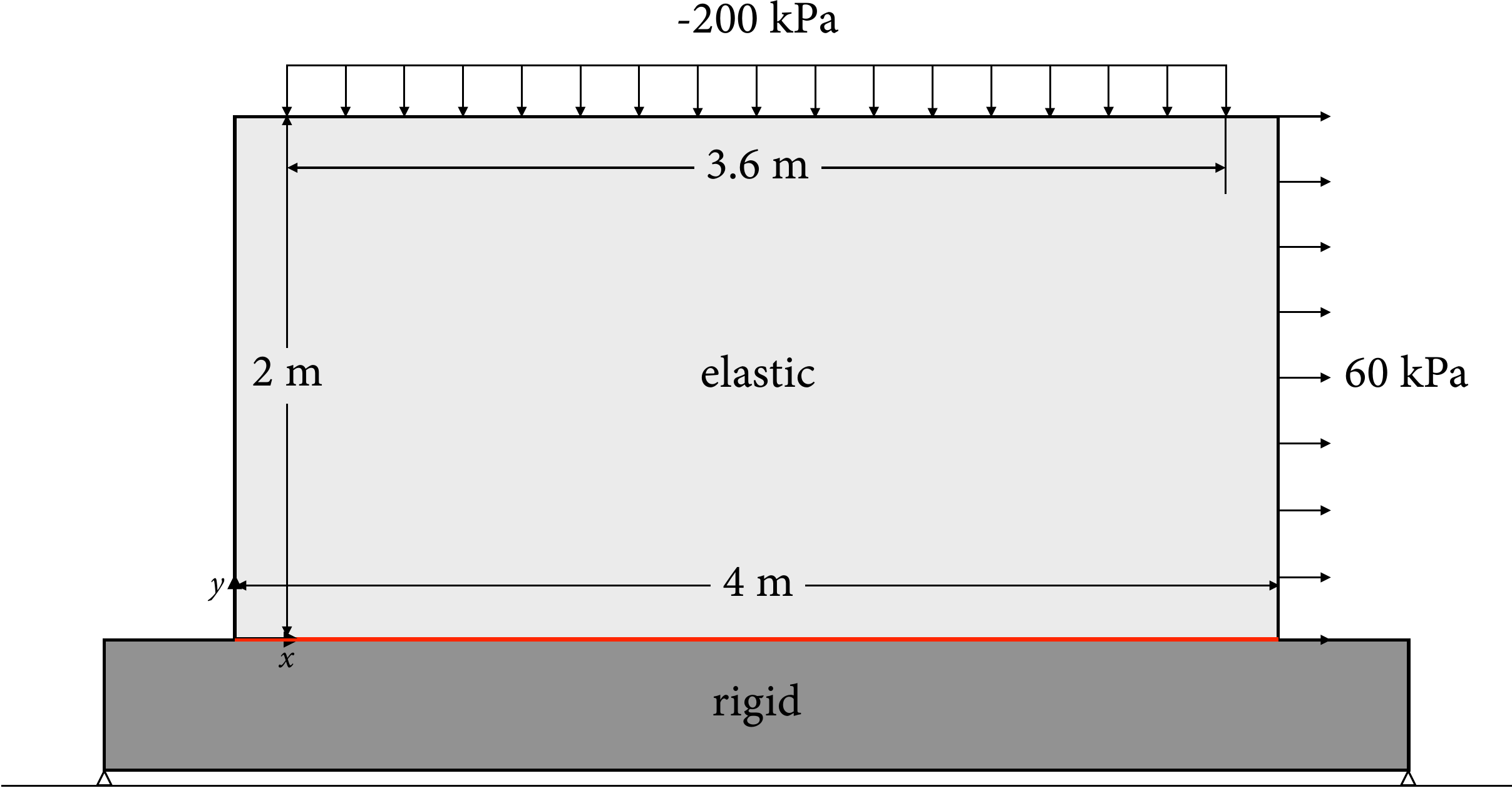}
  \caption{Setup of the sliding of a block problem.}
  \label{fig:sliding-setup}
\end{figure}

In our third example, we simulate a problem whereby stick and slip conditions coexist along an interface.
The problem is an elastic block sliding on a rigid foundation, which was introduced by Oden and Pires~\cite{Oden1984} and later used by other works such as Wriggers \etal~\cite{Wriggers1990} and Simo and Laursen~\cite{Simo1992}.
A slightly modified but essentially the same problem was also presented in Annavarapu \etal~\cite{Annavarapu2014}.
As depicted in Fig.~\ref{fig:sliding-setup}, this problem considers a rectangular elastic block on a rigid foundation
and applies tractions on its top and right boundaries.
The elasticity parameters of the block are $E=1000$ kPa and $\nu=0.3$.
The rigid foundation is approximated with a $10^{9}$ times larger Young's modulus and a zero Poisson's ratio, as done in Annavarapu \etal~\cite{Annavarapu2014}.
Emulating the setup of the original problem, the interface is frictional with $\mu=0.5$ in the 3.6 m-long middle part, but it is frictionless elsewhere.
Under the given condition, the frictional part will mostly be sticked to the foundation but the frictionless parts will slip.
We test two cases of length parameters, $L=0.016$ m and 0.008 m, with meshes of $L/h=8$.
The problem is solved in a single step as in previous works.

\begin{figure}[htbp]
  \centering
  \includegraphics[width=0.95\textwidth]{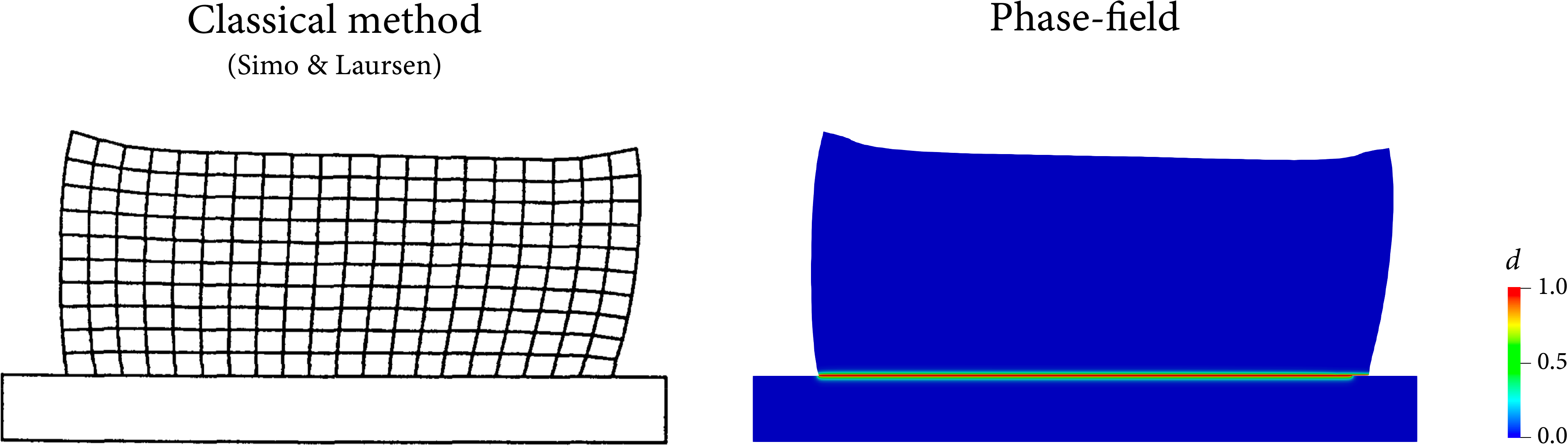}
  \caption{Comparison of deformed geometries obtained by classical finite elements in Simo and Laursen~\cite{Simo1992} and the phase-field method. The phase-field result has been obtained with $L=0.016$ m.}
  \label{fig:sliding-comparison}
\end{figure}

\begin{figure}[htbp]
  \centering
  \includegraphics[width=0.95\textwidth]{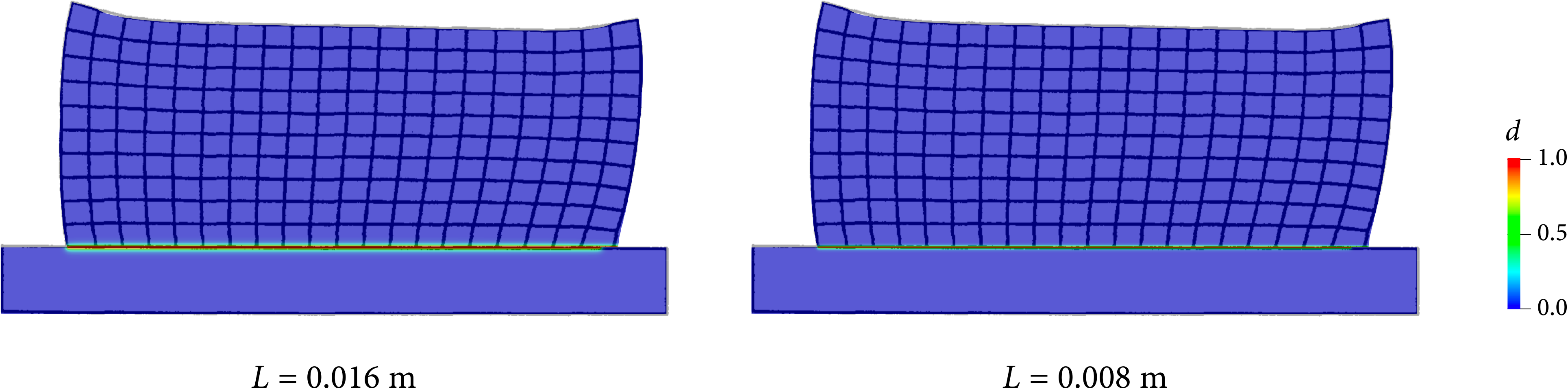}
  \caption{Deformed geometries obtained by the phase-field method with $L=0.016$ m and 0.008 m, superimposed by the deformed geometry in Simo and Laursen~\cite{Simo1992} shown in Fig.~\ref{fig:sliding-comparison}.}
  \label{fig:sliding-comparison-overlapped}
\end{figure}

Figure~\ref{fig:sliding-comparison} compares deformed geometries obtained by our phase-field formulation with numerical result in Simo and Laursen~\cite{Simo1992} obtained by classical finite elements with an augmented Lagrangian method.
It can be seen that the two results are fairly similar and that the interface is partially slipped.
For a more direct comparison, in Fig.~\ref{fig:sliding-comparison-overlapped} the classical result is overlapped to phase-field results obtained with $L=0.016$ m and 0.008 m.
We observe that the two results are matched remarkably well, for both the $L=0.016$ m and 0.008 m cases.
This agreement again demonstrates that the phase-field formulation can correctly identify and reproduce stick and slip behaviors.

\begin{figure}[htbp]
  \centering
  \subfloat[Normal pressures]{\includegraphics[height=0.265\textheight]{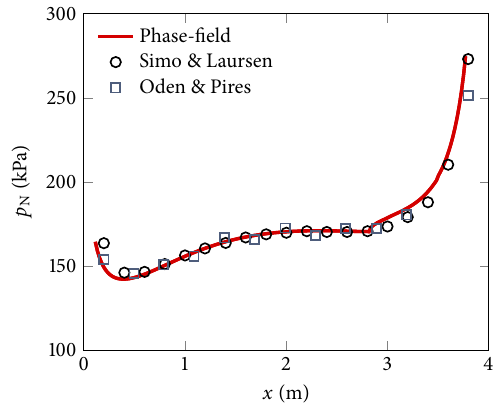}}$\,\,$
  \subfloat[Tangential stresses]{\includegraphics[height=0.265\textheight]{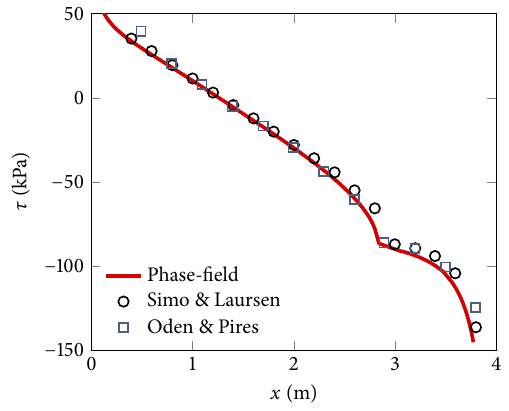}}
  \caption{Contact normal pressures and tangential stresses in comparison with data digitized from Simo and Laursen~\cite{Simo1992} and Oden and Pires~\cite{Oden1984}.}
  \label{fig:sliding-stresses}
\end{figure}

Previous studies have commonly used the contact normal pressures and tangential stresses of this problem to study the performance of contact algorithms.
Here we also use these quantities to fully verify that the phase-field formulation can treat contact constraints well without an algorithm.
Figure~\ref{fig:sliding-stresses} presents the contact normal pressures and tangential stresses at the interface (when $L=0.008$ m), comparing them with data digitized from Simo and Laursen~\cite{Simo1992} and Oden and Pires~\cite{Oden1984}.
One can see that the contact pressures and stresses of the phase-field and classical solutions are also in excellent agreement.
This comparison has verified that the contact stresses and pressures of the phase-field method are accurate even when slip and stick conditions are mixed.
Remarkably, embedded discontinuity methods have sometimes produced oscillatory solutions to the contact tractions of the same problem (see Fig. 8(d) of Annavarapu \etal~\cite{Annavarapu2014} for example).
However, as demonstrated throughout this section, the phase-field method is always free of such oscillation because it reformulates a contact problem as a continuum problem for which no algorithm is necessary for contact constraints.
Therefore, an appealing feature of the phase-field method is that the method can address arbitrary crack geometry as embedded discontinuity methods, but without oscillations in contact pressures and stresses.

\begin{figure}[htbp]
  \centering
  \includegraphics[width=0.55\textwidth]{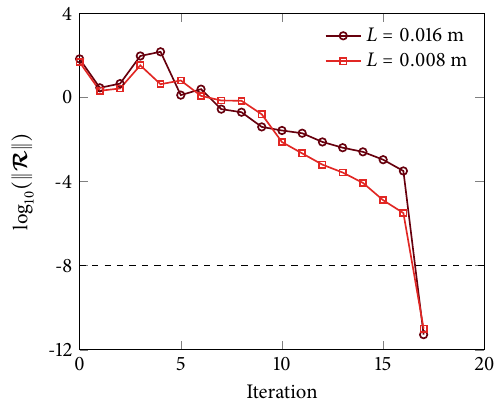}
  \caption{Residual vector norms during Newton iterations in the $L=0.016$ m and $L=0.008$ m cases.}
  \label{fig:sliding-newton}
\end{figure}

Lastly, we plot the Newton convergence behaviors of the $L=0.016$ m and 0.008 m cases in Fig.~\ref{fig:sliding-newton}.
As shown, this problem requires more iterations than prior examples because the interface here involves both stick and slip conditions (the initial contact condition was slip).
Therefore, Newton's method did not converge well initially.
However, once the contact conditions of all points became correctly identified, the residual decreased at a rate for a linear problem.
This behavior agrees well with our observation in the previous example.

\subsection{Propagation of an inclined frictional crack}
Following verification with stationary interface problems,
we simulate propagation of a frictional crack to demonstrate the capability of the phase-field method for modeling growth of a crack with frictional contact.
To this end, we now allow a crack to evolve according to fracture mechanics theory by solving the phase-field equation~\eqref{eq:var-phasefield} in every load step.
The phase-field equation and the momentum balance equation~\eqref{eq:var-mom} are solved sequentially as proposed by Miehe \etal~\cite{Miehe2010a}.
Because this sequential solution method is now fairly standard in the literature, its details are omitted for brevity.
Also, because the phase-field equation has a physical meaning now, $W$ and $G_{c}$ in the phase-field equation should be calculated from physical quantities, rather than being arbitrarily assigned as before.
Considering brittle shear fracture, we regard $W$ as the deviatoric part of strain energy and $G_{c}$ as the mode II fracture energy.
In other words, we have modified a standard phase-field model for brittle fracture to accommodate frictional contact.

\begin{figure}[htbp]
  \centering
  \includegraphics[width=0.4\textwidth]{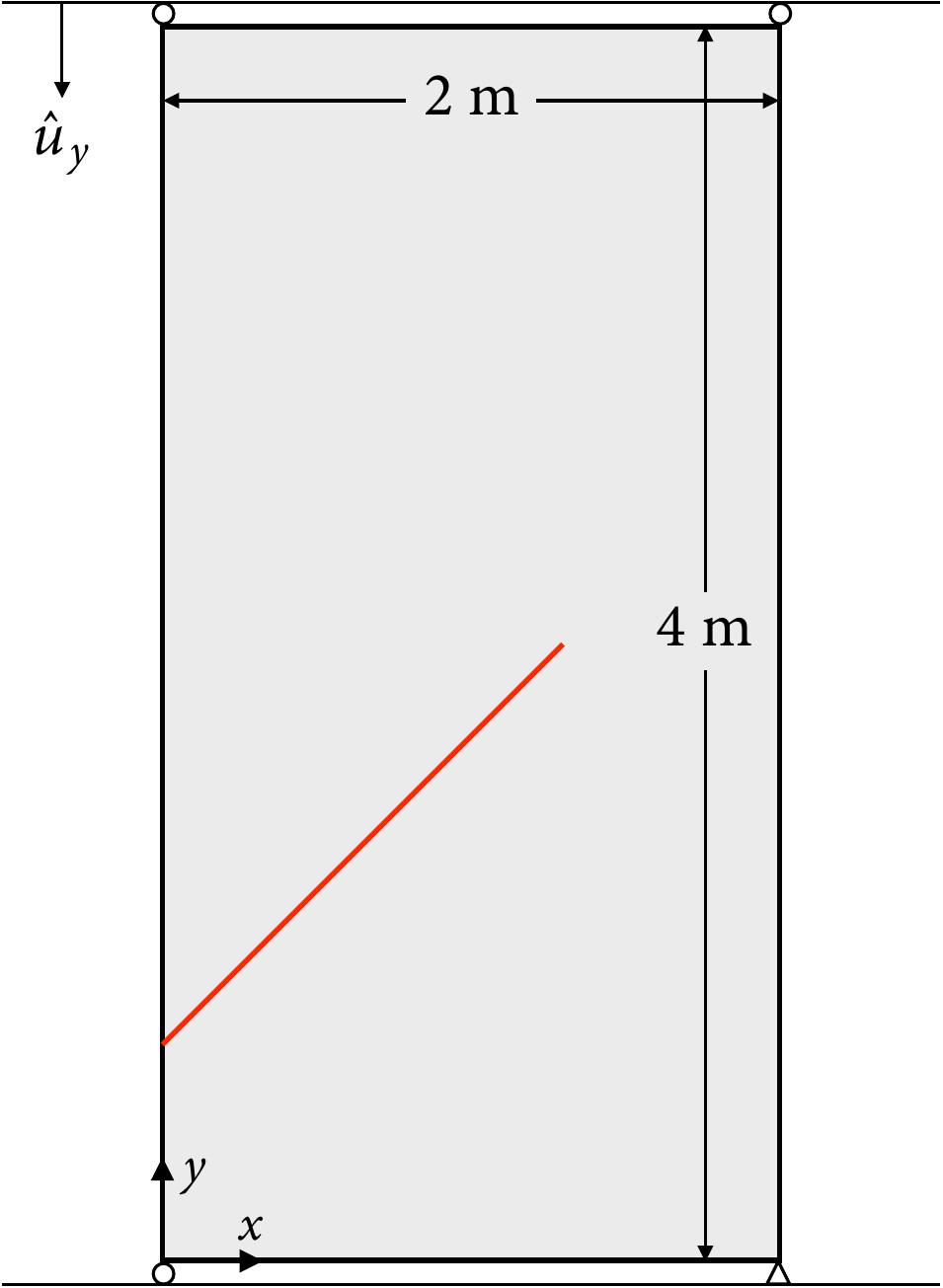}
  \caption{Setup of the propagation of an inclined frictional crack problem.}
  \label{fig:propagation-setup}
\end{figure}

The setup of our particular problem is illustrated in Fig.~\ref{fig:propagation-setup}.
The domain is a 2 m wide and 4 m tall rectangle that possesses a 45$^{\circ}$ inclined crack from coordinates (0.0,0.7) m to (1.3,2.0) m.
The material parameters of the domain are: $E = 10000$ MPa, $\nu = 0.3$, and $G_{c}=50$ kJ/m$^{2}$.
To investigate the effect of friction on this problem, we consider three values of friction coefficients, $\mu=0.01$, 0.10, and 0.30.
For phase-field modeling, we use $L=0.016$ m and locally refine the mesh until $L/h$ reaches 8 along the existing and expected crack path.
Note that the mesh is structured and so the elements are not aligned with the crack direction.
Once the preexisting crack is initialized as before, we vertically compress the domain with a constant displacement rate of $2\times10^{-4}$ m per load step.
Because the contact condition of this problem is rather simple, in most load steps Newton's method converged after a single update.

\begin{figure}[htbp]
  \centering
  \subfloat[Phase field]{\includegraphics[width=0.95\textwidth]{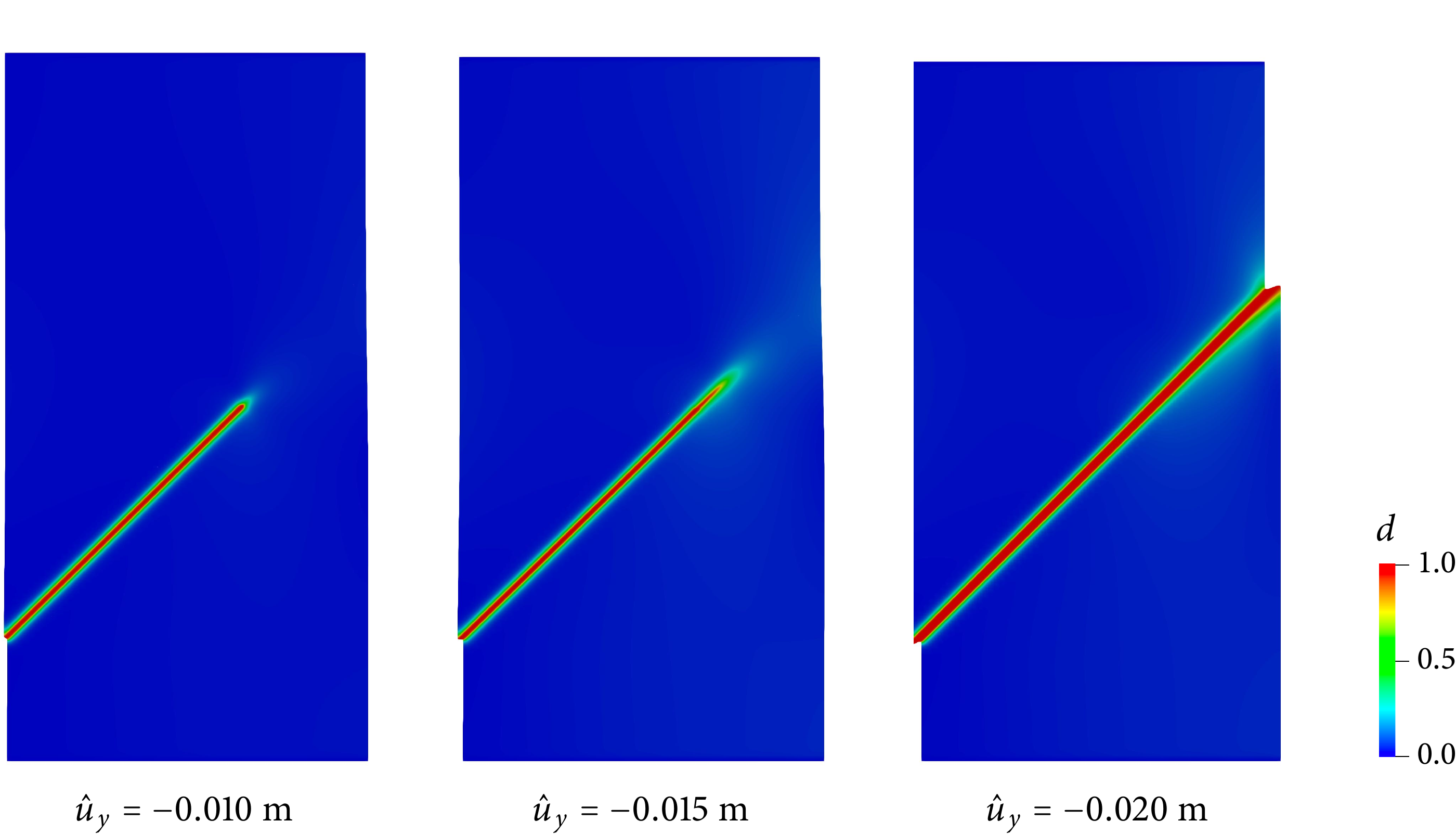}}\\ \vspace{1em}
  \subfloat[$y$-displacement (normalized)]{\includegraphics[width=0.95\textwidth]{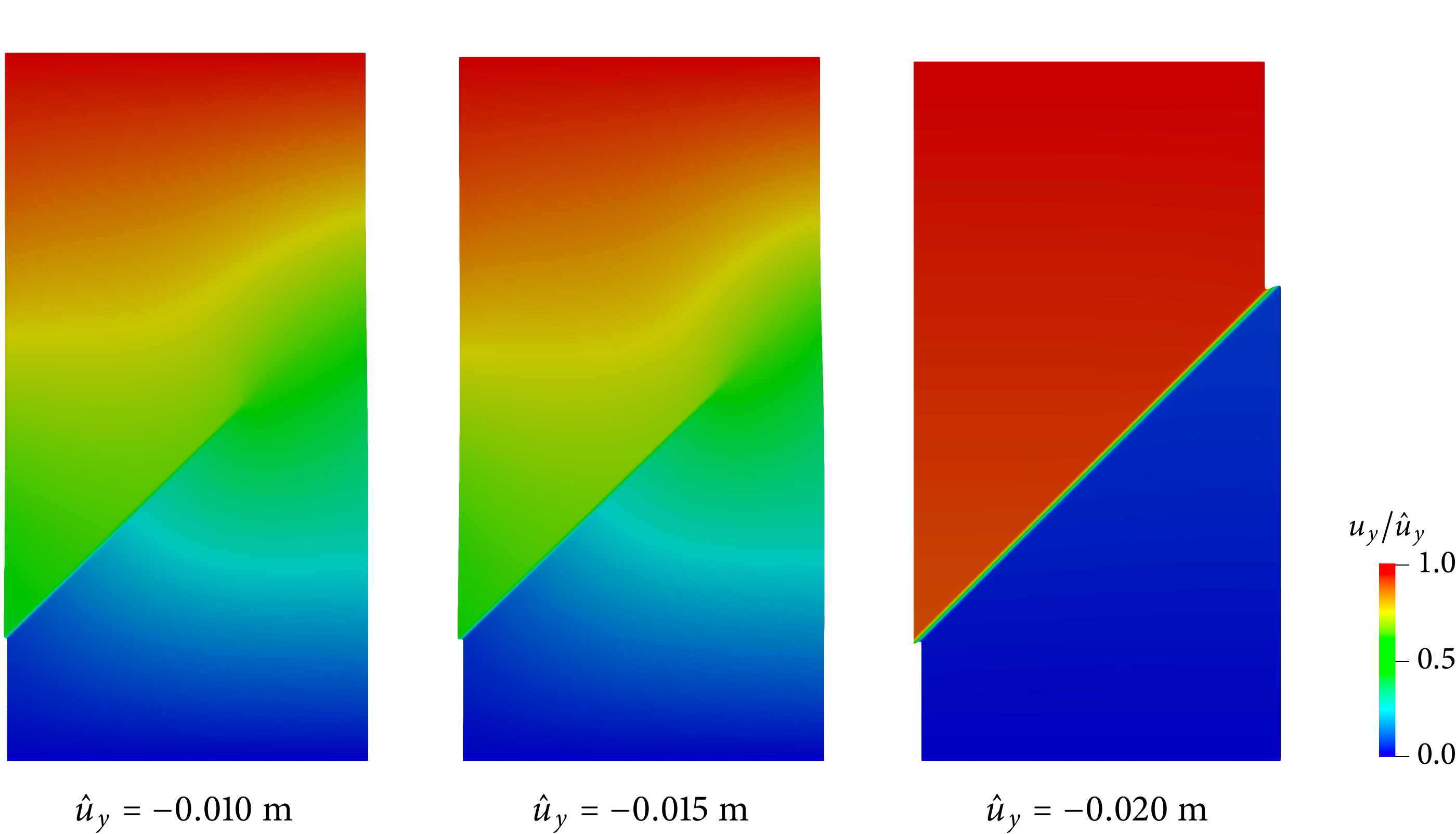}}
  \caption{Simulation results of the crack propagation problem with $\mu=0.30$. Displacement is scaled by a factor of 5.}
  \label{fig:propagation-phase-field}
\end{figure}

Figure~\ref{fig:propagation-phase-field} shows how the phase-field variable and the vertical displacement field evolve during the course of loading when $\mu=0.3$.
As shown, the phase-field model well simulates propagation of the preexisting crack along the 45$^{\circ}$ direction until it reaches the upper right size of the domain.
During the propagation stage, the displacement field is discontinuous across the crack but still continuous through the non-fractured region.
After the crack has fully developed, however, the upper and lower parts of the domain are completely disconnected, and the upper part slips along the crack.
Note that this post-fracture process is essentially the same as stationary interface problems simulated earlier in this section.
We also note that other friction coefficient cases show qualitatively identical responses in terms of the crack path and the displacement pattern.

\begin{figure}[htbp]
  \centering
  \includegraphics[width=0.55\textwidth]{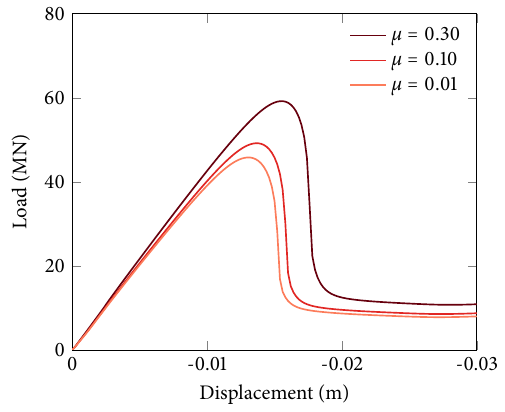}
  \caption{Load--displacement curves of the three friction coefficient cases.}
  \label{fig:propagation-load-disp}
\end{figure}

In Fig.~\ref{fig:propagation-load-disp} we plot the load--displacement curves of the three friction coefficient cases.
As expected, the peak load and displacement increase with the friction efficient.
We also see that the two cases show more or less the same pattern, in which the material fails in a brittle manner and then exhibits a residual strength.
The residual strength also increases with the friction coefficient, which evinces the contribution from the frictional resistance along the crack.

Before closing this section, we would like to demonstrate the critical role of contact treatment in phase-field modeling of crack propagation under compression.
For this purpose, we simulate the same problem with the model of Amor \etal~\cite{Amor2009}, whereby the contact condition is treated by the volumetric--deviatoric decomposition of the stress tensor.
This stress decomposition scheme is the only difference from our phase-field formulation used above.
We note that when our formulation attempted to simulate this problem without friction ($\mu=0$), it did not converge from the very first load step because all nodes along the preexisting crack slip immediately.
However, although the volumetric--deviatoric decomposition assumes frictionless contact, it can still simulate this problem until the crack fully develops.
This indicates that the inexact contact treatment of the volumetric--deviatoric decomposition provides non-physical frictional resistance along the interface.
On a related note, we have also found that the volumetric--deviatoric decomposition is unable to distinguish stick and slip conditions correctly for the second example of this section.

\begin{figure}[htbp]
  \centering
  \subfloat[Phase field\vspace{1em}]{\includegraphics[width=0.8\textwidth]{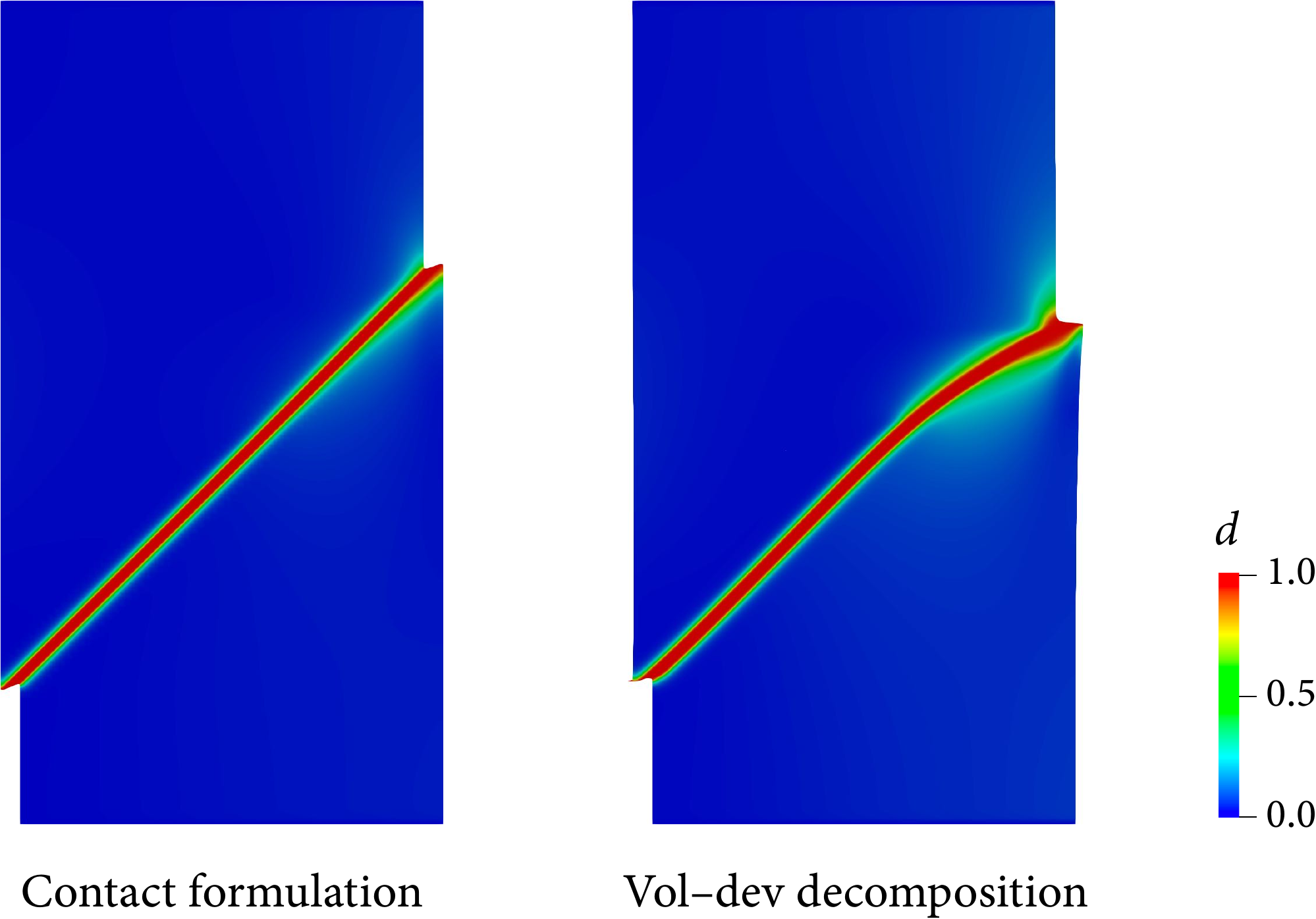}}\\ \vspace{1em}
  \subfloat[$y$-displacement (normalized)\vspace{1em}]{\includegraphics[width=0.8\textwidth]{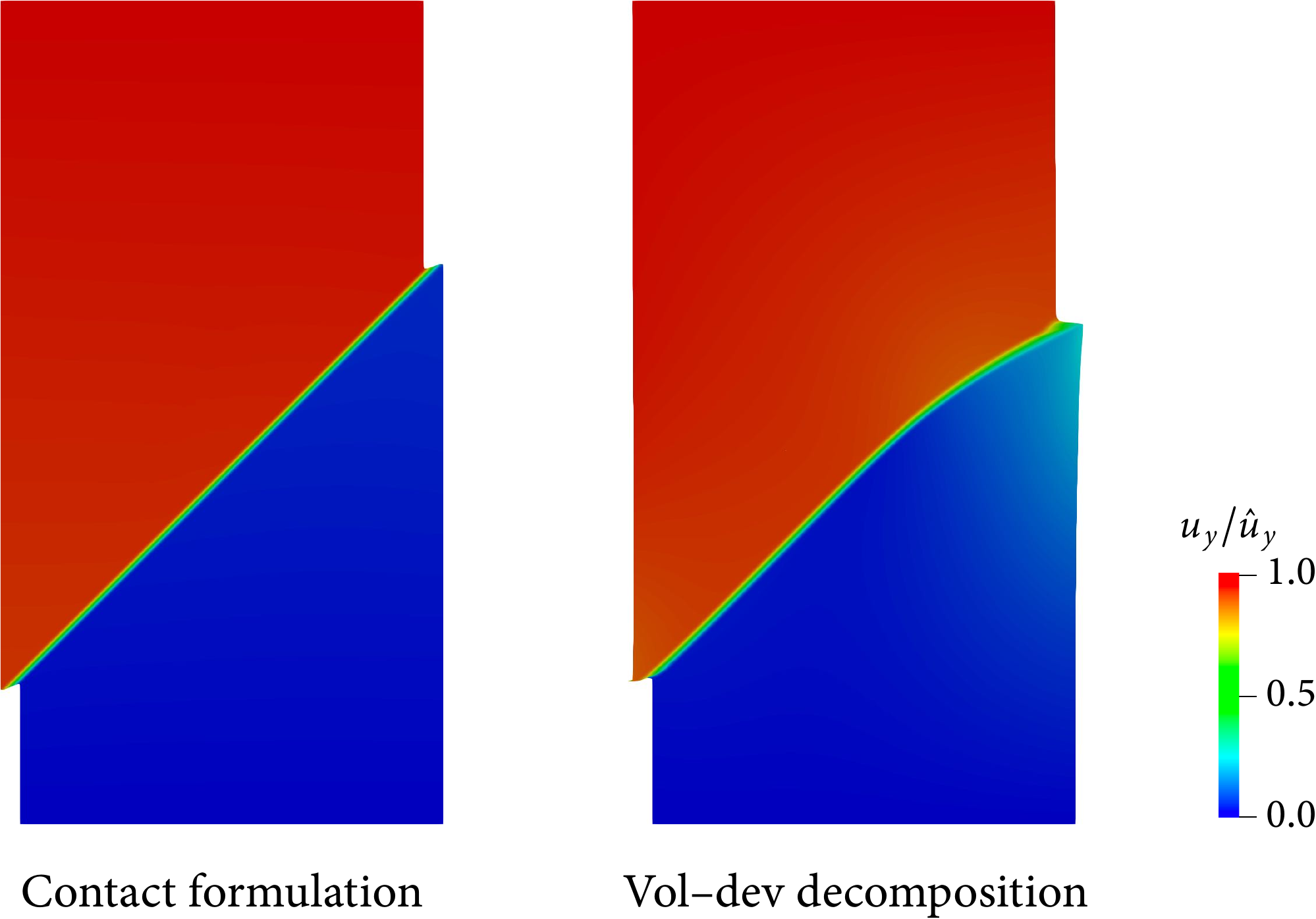}}
  \caption{Comparison of simulation results at $\hat{u}_{y}=-0.020$ m between those obtained by phase-field models with the contact formulation proposed in this work (with $\mu=0.01$) and the volumetric--deviatoric decomposition proposed by Amor \etal~\cite{Amor2009}. Displacement is scaled by a factor of 5.}
  \label{fig:propagation-comparison}
\end{figure}

Figure~\ref{fig:propagation-comparison} compares simulation results from the two phase-field formulations when $\hat{u}_{y}=-0.020$ m.
To minimize the influence of friction on the comparison, the $\mu=0.01$ case is shown in this figure.
One can see that when the volumetric--deviatoric decomposition is used for this problem, the crack path becomes kinked, giving rise to rather unrealistic deformation responses.
This difference demonstrates that inappropriate estimation of contact stresses can impact the crack driving force to the extent that alters the crack path direction.
Therefore, it can be concluded that accurate treatment of contact condition is critical to the application of phase-field modeling to compression-induced fracture propagation, which is a classic problem in geomechanics~\cite{Hoek1965,Ingraffea1980,Nemat-Nasser1982,Horii1985}.

\section{Closure}
\label{sec:closure}
A phase-field method has been proposed for modeling cracks with frictional contact.
Built on standard approaches in phase-field modeling of fracture, the proposed method calculates the stress tensor in a regularized interface region by identifying the contact condition in an interface-oriented coordinate system.
By doing so, the phase-field method accommodates contact behavior in the interface direction while imposing no-penetration constraints in other directions.
Using benchmark examples in the literature, we have verified that the proposed method can provide numerical solutions very close to those obtained by a discrete method, showing little sensitivity to the length parameter for phase-field regularization.
Moreover, by allowing the crack to evolve according to brittle fracture theory, we have demonstrated that the proposed phase-field method can also simulate propagation of frictional cracks.

The proposed phase-field method has two key features that make it an appealing alternative to standard discrete methods.
First, it can model a crack passing through the interior of elements without an explicit representation of geometry or enrichment of basis functions.
Second, it does not require a sophisticated algorithm for imposing contact constraints on crack surfaces.
Thanks to these two features, the phase-field method can be implemented far more easily than most of existing methods for frictional cracks.
\revised{Furthermore, much like phase-field models of opening fractures, the proposed method may be well applied to more complex settings such as finite deformations, bulk plasticity, and coupled multiphysics.}
Therefore, it is believed that the phase-field method can be an attractive option even for modeling frictional interfaces that are stationary, let alone those that evolve.

\section*{Acknowledgments}
This work was supported by the Research Grants Council of Hong Kong (Project 27205918).
The first author also acknowledges financial support from the Hong Kong PhD Fellowship.

\bibliography{references}

\end{document}